\documentclass[final]{acmsiggraph}

\usepackage[ruled]{algorithm2e} 

\usepackage{simplewick}
\usepackage{tikz}
\usepackage{amsfonts}
\usepackage{amsmath}
\usepackage{amsthm}
\usepackage{amssymb}
\usepackage{float}
\usepackage{subfig}
\usepackage{tabularx}
\usepackage{epstopdf}
\usepackage{centernot}
\usepackage[utf8]{inputenc}
\usepackage{cleveref}
\usepackage{stmaryrd}
\usepackage{soul}
\usepackage{wrapfig}
\usepackage{lipsum}
\usepackage{xcolor,colortbl}
\usepackage{amssymb}
\usepackage{pifont}
\usepackage{booktabs}
\usepackage{colortbl}
\usepackage{xfrac}
\usepackage{caption}
\usepackage{wrapfig}
\usepackage{lipsum}
\usepackage{algorithmic}
\usepackage{enumitem}

\usepackage{placeins}

\newcommand{\vc}[1]{\mbox{\textbf{{$\mathsf #1$}}}}

\newcommand{\R}[0]{\mathbb{R}}

\definecolor{redcolor}{rgb}{0.8,0,0}
\definecolor{bluecolor}{rgb}{0.0,0.1,0.6}
\definecolor{orangecolor}{rgb}{0.9.,0.5.,0.1}
\definecolor{greencolor}{rgb}{0.5,0.7,0.5}
\definecolor{browncolor}{rgb}{0.5,0.2,0.2}
\definecolor{greycolor}{rgb}{0.6,0.6,0.6}

\def\argmin{\mathop{\rm argmin}}
\def\min{\mathop{\rm min}}

\newcommand{\red}[1]{{\color{redcolor}{#1}}}
\newcommand{\orange}[1]{{\color{orangecolor}{#1}}}
\newcommand{\blue}[1]{{\color{bluecolor}{#1}}}

\newcommand{\bfi}[1]{\textit{ \textbf{#1}}}

\title{Blended Cured Quasi-Newton for Geometry Optimization}

\author{
 Yufeng Zhu \\{\footnotesize University of British Columbia \& Adobe Research}
  \and Robert Bridson \\{\footnotesize Autodesk}
  \and Danny M. Kaufman \\{\footnotesize Adobe Research}
}

\pdfauthor{}

\begin{document}

\teaser{
   \includegraphics[width=\textwidth]{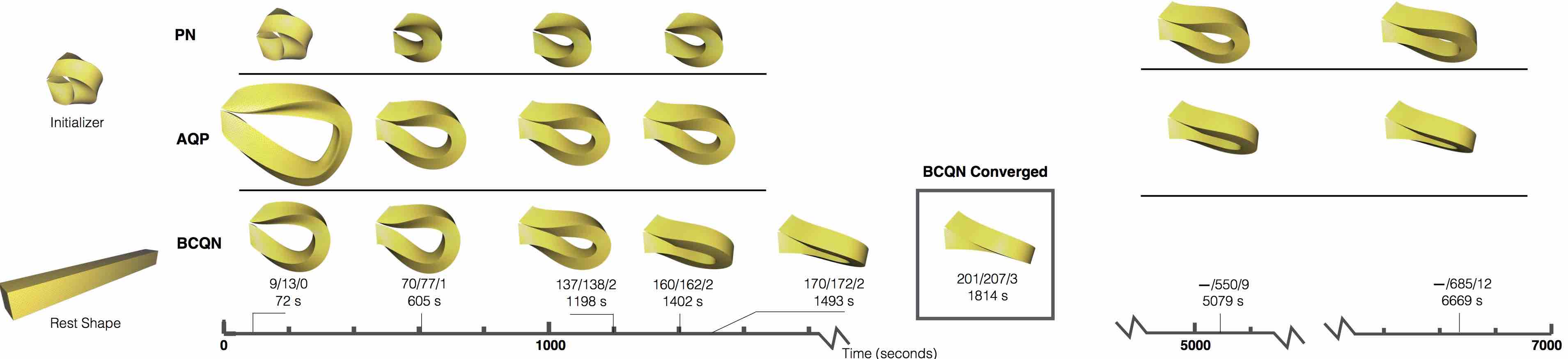}
 	\caption{\bfi{Twisting.} A stress-test 3D deformation problem. {\bf Left:} we initialize a 1.5M tetrahedra mesh bar with a straight rest shape into a tightly twisted coil, constraining
both ends to stay fixed.  {\bf Right:} minimizing the ISO deformation energy to find a constrained equilibrium with (top to bottom) Projected Newton (PN), Accelerated Quadratic Proxy (AQP) and our BCQN method, we show intermediate shapes at reported wall-clock time (seconds) and iteration counts at those times (BCQN/AQP/PN). BCQN converges at 30 minutes, while AQP and PN continue to optimize.}
 	\label{fig:teaser}
}

\maketitle

\begin{abstract}
Optimizing deformation energies over a mesh, in two or three
dimensions, is a common and critical problem in physical simulation and geometry processing. 
We present three new improvements to the
state of the art: a barrier-aware line-search filter that cures blocked descent steps due to  
element barrier terms and so enables rapid progress;
an energy proxy model that adaptively blends the Sobolev (inverse-Laplacian-processed) gradient and L-BFGS descent to gain the advantages of both,  
while avoiding L-BFGS's current limitations in geometry optimization tasks;
and a characteristic gradient norm providing a
robust and largely mesh- and energy-independent convergence criterion
that avoids wrongful termination when algorithms temporarily slow
their progress. Together these improvements form the basis for
Blended Cured Quasi-Newton (BCQN), a new geometry optimization
algorithm. Over a wide range of problems
over all scales we show that BCQN is generally the fastest and most
robust method available, making some previously intractable problems
practical while offering up to an order of magnitude improvement
in others.
\end{abstract}

\section{Introduction}
Many fundamental \emph{physical} and \emph{geometric} modeling tasks
reduce to minimizing nonlinear measures of deformation over meshes.
Simulating elastic bodies, parametrization, deformation, shape
interpolation, deformable inverse kinematics, and animation all
require robust, efficient, and easily automated \emph{geometry
optimization}.
By \emph{robust} we mean the algorithm should solve
every reasonable problem to any accuracy given commensurate time,
and only reports success when the accuracy has truly been achieved.
By \emph{efficient} we mean rapid convergence in wall-clock time,
even if that may mean more (but cheaper) iterations.
By \emph{automated} we mean the user needn't adjust algorithm
parameters or tolerances at all to get good results
when going between different problems. With these three attributes,
a geometry optimization algorithm can be reliably used in production
software.

We propose a new algorithm, Blended Cured Quasi-Newton (BCQN),
with three core contributions based on analysis of where prior methods
faced difficulties:
\begin{itemize}
\item an \textbf{adaptively blended quadratic energy proxy} for the
deformation energy combining the Sobolev gradient and a quasi-Newton
secant approximation, allowing a low cost per iterate with second-order
acceleration but avoiding secant artifacts where the Laplacian is
more robust;
\item a \textbf{barrier-aware filter} 
on search directions, that gains larger step sizes and so improved convergence progress in line search for the common case of iterates where individual elements degenerate towards collapse.
\item a \textbf{characteristic gradient norm convergence criterion},
which is immune to terminating prematurely due to algorithm stagnation,
and is consistent across mesh sizes, scales, and choice of energy
so per-problem adjustment is unnecessary.
\end{itemize}
Over a wide range of test cases we show that BCQN makes the solution
of some previously intractable problems practical, offers up to an
order of magnitude speed-up in other cases, and in all cases
investigated so far either improves on or closely matches the
performance of the best-in-class optimizers available. We claim
BCQN achieves our goals for production software.

\section{Problem Statement and Overview}

The geometry optimization problem we face is solving
\begin{equation}
    x^* = \argmin_{x\in \R^{dn}} E(x),
\end{equation}
for $n$ vertex locations in $d$-dimensional space stored in vector $x$,
where the energy $E(x)$ is a measure of the deformation, and $x$
is subject to boundary conditions.\footnote{We restrict our attention to
constraining a subset of vertex positions to given values, i.e.\ Dirichlet conditions,
for simplicity.} The energy is expressed as a sum over elements $t$ in a triangulation $T$
(triangles or tetrahedra depending on dimension),
\begin{equation}
\label{eq:obj}
E(x) = \sum_{t \in T} a_t W\big( F_t(x) \big),
\end{equation}
where $a_t>0$ is the area or volume of the rest shape of element $t$, $W$ is an energy
density function taking the deformation gradient as its argument, and $F_t$ computes the
deformation gradient for element $t$.
This problem may be given as is, or may be the result of a discretization of
a continuum problem with linear finite elements for example.

\subsection{Iterative solvers for nonlinear minimization}

Solution methods for the above generally apply an algorithmic strategy of iterated
approximation and stepping~\cite{Bertsekas:2016:NOP}, built
from three primary ingredients: an energy approximation, a
line search, and a termination criteria.\footnote{Alternatively, trust-region
methods are available, though not considered in the current work nor as
popular within the field.} \\

\bfi{Energy Approximation} At the current iterate $x_i$ we form a
local quadratic approximation of the energy, or \emph{proxy}:
\begin{align}
\label{eq:quad_approx}
E_i(x) = E(x_i) +   (x - x_i)^T \nabla E(x_i)  + \tfrac{1}{2}  (x - x_i) ^T H_i (x - x_i)
\end{align} where $H_i$ is a symmetric matrix.
Near the solution, if $H_i$ accurately approximates the Hessian we can achieve
fast convergence optimizing this proxy, but it is also critical that
it be stable --- symmetric positive definite (SPD) --- to ensure the proxy
optimization is well-posed everywhere; we also want $H_i$ to be cheap to solve with,
preferring sparser matrices and ideally not having to refactor at each iteration. \\

\bfi{Line Search} Quadratic models allow us to apply linear solvers
to find stationary points $x_i^* = \argmin_x E_i(x)$ of the local
energy approximation. A step
\begin{align}
\label{eq:descent_step_solve}
p_i = x_i^* - x_i = -H_i ^{-1} \nabla E(x_i) 
\end{align}
towards this stationary point then forms a direction for probable energy descent.
However, quadratic models are only locally accurate for nonlinear energies in general,
thus line-search is used to find an improved length $\alpha_i>0$ along $p_i$ to get a new iterate 
\begin{align}
\label{eq:vanilla_step}
x_{i+1} \leftarrow x_i + \alpha_i p_i,
\end{align}
for adequate decrease in nonlinear energy $E$. Of particular concern for
the geometric problems we face is energies which blow up to infinity for
degenerate (flattened) elements: in a given step, the elements where this
may come close to happening rapidly depart from the proxy, and the step size $\alpha_i$ may
have to be very small indeed, see Figure\ \ref{fig:blocked_line_search}, impeding progress globally. \\

\bfi{Termination} Iteration continues until we are able to stop with
a ``good enough'' solution -- but this requires a precise computational
definition. Typically we monitor some quantity which approaches zero
if and \emph{only if} the iterates are approaching a stationary point.
The standard in unconstrained optimization is to check the norm of the
gradient of the energy, which is zero only at a stationary point and
otherwise positive; however, the raw gradient norm depends on the mesh
size, scaling, and choice of energy, which makes finding an appropriate
tolerance to compare against highly problem-dependent and difficult
to automate. \\

\section{Related Work}

\subsection{Energies and Applications}

A wide range of physical simulation and geometry processing
computations are cast as \emph{variational} tasks to minimize
measures of distortion over domains.

To simulate elastic solids with large deformations,
we typically need to minimize hyper-elastic potentials formed
by integrating strain-energy densities over the body. These
material models date back to Mooney~\shortcite{Mooney:1940:ATO} and
Rivlin~\shortcite{Rivlin:1948:SAO}.  Their Mooney-Rivlin and
Neo-Hookean materials, and many subsequent hyperelastic materials,
e.g.~St.~Venant-Kirchoff, Ogden,  Fung~\cite{Bonet:1998:ASO}, are
constructed from empirical observation and
analysis of deforming real-world materials. Unfortunately, all but
a few of these energy densities are nonconvex. This makes their
minimization highly challenging. Constants in these models are
determined by experiment for scientific computing
applications~\cite{Ogden:1972:LDI}, or alternately are directly set
by users in other cases~\cite{Xu:2015:NMD}, e.g., to meet artistic
needs.

In geometry processing a diverse range of energies have 
been proposed to minimize various mapping distortions, 
generally focused on minimizing either measures of
isometric~\cite{Sorkine:2007:ARA,Chao:2010:ASG,Smith:2015:BPW,Aigerman:2015:Seamless,Liu:2008:ALG}
or
conformal~\cite{Hormann:2000:MIPS,Levy:2002:LSC,Desbrun:2002:IPO,Benchen:2008:CFB,Mullen:2008:SCP,Weber:2012:CEQ}
distortion. While some of these energies do not prohibit
inversion~\cite{Sorkine:2007:ARA,Chao:2010:ASG,Levy:2002:LSC,Desbrun:2002:IPO}
many others have been explicitly constructed with nonconvex terms
that guarantee preservation of local
injectivity~\cite{Hormann:2000:MIPS,Aigerman:2015:Seamless,Smith:2015:BPW}.
Other authors have also added constraints to strictly bound distortion, for example,
but we restrict attention to unconstrained minimization --- but note constrained
optimization often relies on unconstrained algorithms as an inner kernel.

Our goal here is to provide a tool to minimize arbitrary energy
density functions as-is. We take as input energy functions provided
by the user, irrespective of whether these energies are custom-constructed
for geometry tasks, physical energies extracted from experiment,
or energies hand-crafted by artists. Our work focuses on the better
optimization of the important \emph{nonconvex} energies whose
minimization remains the primary challenging bottleneck in many
modern geometry and simulation pipelines.  In the following sections,
to evaluate and compare algorithms, we consider a range of challenging
nonconvex deformation energies currently critical in physical
simulation and geometry processing: Mooney-Rivlin
{\bf(MR)}~\cite{Bower:2009:AMO}, Neo-Hookean
{\bf(NH)}~\cite{Bower:2009:AMO},  Symmetric Dirichlet
{\bf(ISO)}~\cite{Smith:2015:BPW}, Conformal Distortion
{\bf(CONF)}~\cite{Aigerman:2015:Seamless}, and Most-Isometric
Parameterizations {\bf(MIPS)}~\cite{Hormann:2000:MIPS}.

\subsection{Energy Approximations} 

Broadly, existing models for the local energy approximation in (\ref{eq:quad_approx}) fall into four rough categories
that vary in the construction of the \emph{proxy}\footnote{Names and notations for $H_i$ vary across the literature
depending on method and application. For consistency, here, across all
methods we will refer to $H_i$ as the \emph{proxy} matrix --- inclusive of cases where it is the actual Hessian
or direct modification thereof.} matrix $H_i$.
\emph{Newton-type} methods exploit expensive
second-order derivative information;
\emph{first-order} methods use only first derivatives and
apply lightweight fixed proxies;
\emph{quasi-Newton} methods iteratively update proxies to approximate
second derivatives using just differences in gradients;
\emph{Geometric Approximation} methods use
more domain knowledge to directly construct proxies which relate to 
key aspects of the energy, resembling Newton-type methods but
not necessarily taking second derivatives.

\bfi{Newton-type} methods generally can achieve the most rapid convergence
but require the costly assembly, factorization and backsolve of new
linear systems per step.  At each iterate Newton's method
uses the energy Hessian, $\nabla^2 E(x_i)$, to form a proxy matrix.
This works well for convex energies like ARAP\ \cite{Chao:2010:ASG},
but requires modification for nonconvex energies\ \cite{Nocedal:2006:Book}
to ensure that the proxy is at least positive semi-definite (PSD).
Composite Majorization (CM), a tight convex majorizer, was recently
proposed as an analytic PSD approximation of the Hessian\
\cite{Shtengel:2017:GOV}. The CM proxy is efficient to assemble but
is limited to two-dimensional problems and just a trio of energies:
ISO, NH and symmetric ARAP.  More general-purpose solutions include
adding small multiples of the identity, and projection of the Hessian
to the PSD cone but these generally damp convergence too much\
\cite{Liu:2016:TRT,Shtengel:2017:GOV,Nocedal:2006:Book}.  More
effective is the Projected Newton (PN) method that projects per-element
Hessians to the PSD cone prior to assembly\ \cite{Teran:2005:RQF}.
Both CM and PN generally converge rapidly in the nonconvex setting
with CM often outperforming PN in the subset of 2D cases where CM
can be applied\ \cite{Shtengel:2017:GOV}, while PN is more general
purpose for 3D and 2D problems.  Both PN and CM have identical
per-element stencils and so identical proxy structures. Despite low
iteration counts they both scale prohibitively due to per-iteration
cost and storage as we attempt increasingly large optimization
problems.

\bfi{First-order} methods build descent steps by preconditioning the gradient with a fixed proxy matrix. These proxies are generally inexpensive to solve and sparse so that cost and storage remain tractable as we scale to larger systems. However, they often suffer from slower convergence as we lack higher-order information.
Direct gradient descent, $H_i \leftarrow Id$, and Jacobi-preconditioned gradient descent, $H_i \leftarrow diag(\nabla^2 E(x_i)\big)$ offer attractive opportunities for parallelization~\cite{Wang:2016:DMF,Fu:2015:CLI} but suffer from especially slow convergence due to poor scaling.
The Laplacian matrix, $L$, forms an excellent preconditioner, that both smooths and rescales the gradient~\cite{Neuberger:1985:SDA,Martin:2013:ENL,Kovalsky:2016:AQP}. Unlike the Hessian, the Laplacian is a constant PSD proxy that can be pre-factorized once and backsolved separately per-coordinate. Iterating descent with $H_i \leftarrow L$, is the Sobolev-preconditioned gradient descent (SGD) method. SGD was first introduced, to our knowledge, by Neuberger~\shortcite{Neuberger:1985:SDA}, but has since been rediscovered in graphics as the local-global method for minimizing ARAP~\cite{Sorkine:2007:ARA}. As noted by Kovalsky et al.\ \shortcite{Kovalsky:2016:AQP} Local-global for ARAP is exactly SGD.
More recently Kovalsky et al.~\shortcite{Kovalsky:2016:AQP} developed the highly effective Accelerated Quadratic Proxy (AQP) method by adding a Nesterov-like acceleration~\cite{Nesterov:1983:AMO} step to SGD. This improves AQP's convergence over SGD. However, as this acceleration is applied after line search, steps do not guarantee energy decrease and can even contain collapsed or inverted elements --- preventing further progress. More generally, the Laplacian is
constant and so ignores valuable local curvature information ---
we see this issue in a number of examples in Section~\ref{sec:results}
where AQP stagnates and is unable to converge. Curvature can make
the critical difference to enable progress.

\bfi{Quasi-Newton} methods lie in between these two extremes. They
successively, per descent iterate, update approximations of the
system Hessian using a variety of strategies.  Quasi-Newton methods
employing sequential gradients to updates proxies, i.e.  L-BFGS and
variants, have traditionally been highly successful in scaling up
to large systems~\cite{Bertsekas:2016:NOP}. Their updates can be
performed in a compute and memory efficient manner and
can guarantee the proxy is SPD even where the exact Hessian is not.
While not fully second-order, they achieve superlinear convergence, regaining
much of the advantage of Newton-type methods. L-BFGS convergence
can be improved with the choice of initializer. Initializing with the
diagonal of the Hessian\ \cite{Nocedal:2006:Book}, application-specific
structure\ \cite{Jiang:2004:APL} or even the Laplacian\ \cite{Liu:2016:TRT}
can help. However, so far, for geometry optimization problems,
L-BFGS has consistently and surprisingly failed to perform
competitively~\cite{Kovalsky:2016:AQP,Rabinovich:2016:SLI}
\emph{irrespective} of choice of initializer. Nocedal and Wright point
out that the secant approximation can implicitly create a \emph{dense}
proxy, unlike the sparse true Hessian, directly and incorrectly
coupling distant vertices. This is visible as swelling artifacts
for intermediate iterations in Figure \ref{fig:quadratic_compare}.

\begin{figure}[h!]
\centering
\includegraphics[width=0.9\linewidth]{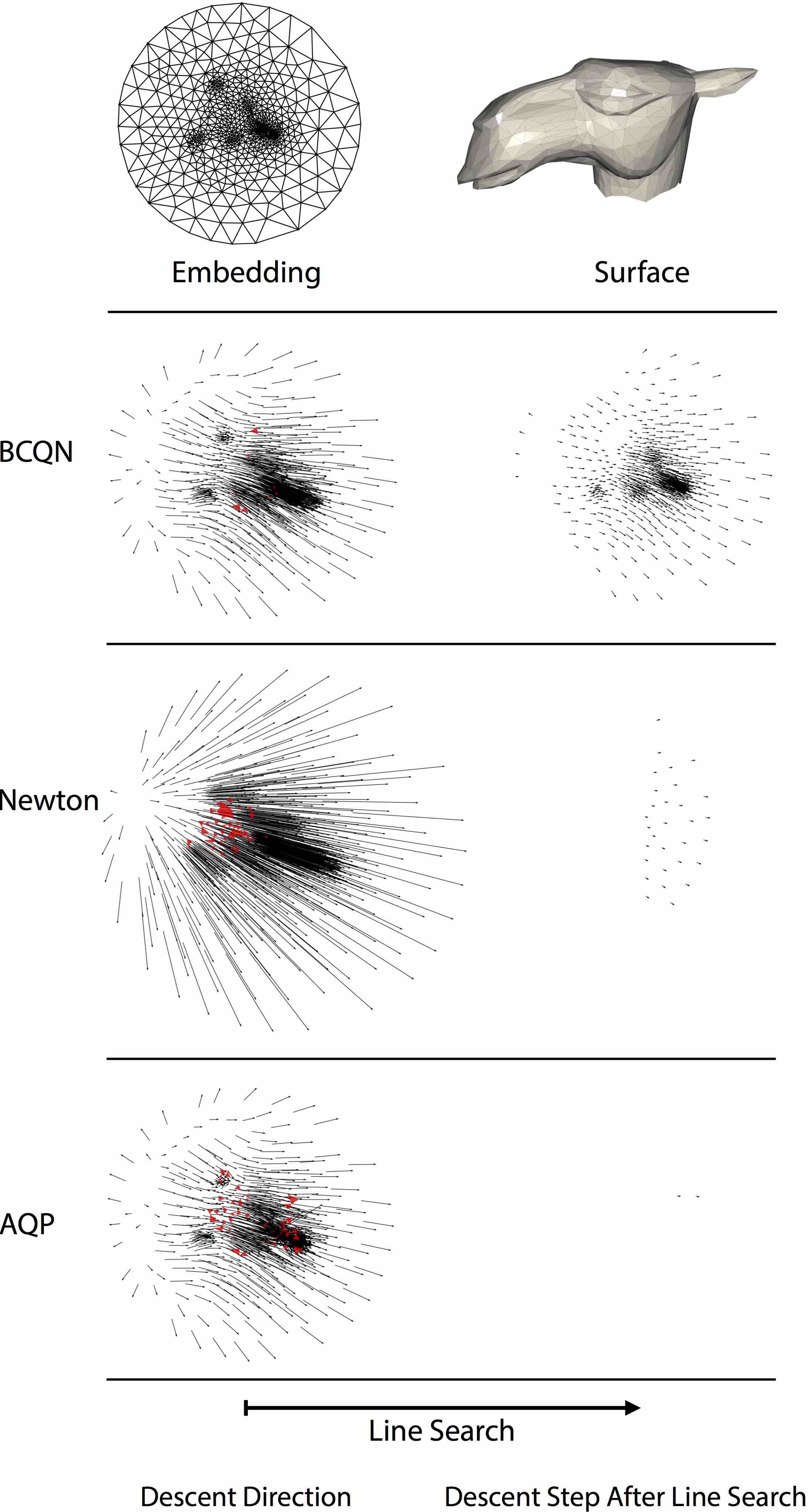}
\caption{\bfi{Line-search blocking.} Barrier terms in nonconvex energies (here we use ISO) of the form $1/g(\sigma)$ can severely restrict step sizes in line searches even when expensive, high-quality methods such as Newton-type methods are applied. {\bf Left column:} descent-direction vector fields, per vertex, in a descent step generated by BCQN, PN and AQP with potential blocking triangles rendered in red. {\bf Right, bottom rows:} after line-search, close to collapsing elements have restricted the global step size for AQP and PN to effectively block progress. {\bf Right, top row:} BCQN's barrier-aware line-search filtering enables progress with significant descent directions.}
\label{fig:blocked_line_search}
\end{figure}

\bfi{Geometric Approximation} methods specifically for geometry optimization
have also been developed recently: SLIM~\cite{Rabinovich:2016:SLI}
and the AKAP preconditioner~\cite{Claici:2017:IAP}. SLIM extends
the local-global strategy to a wide range of distortion energies
while AKAP applies an approximate Killing Vector Field operator as
the proxy matrix. Both require re-assembly and factorization of
their proxies for each iterate. SLIM and AKAP convergence are
generally well improved over SGD and AQP~\cite{Rabinovich:2016:SLI,Claici:2017:IAP}. However, they do not match the convergence quality
of the second-order, Newton-type methods, CM and
PN~\cite{Shtengel:2017:GOV}. SLIM falls well short of both CM and
PN~\cite{Shtengel:2017:GOV}.  AKAP is more competitive than SLIM
but remains generally slower to converge than PN in our testing,
and is much slower than CM.  At the same time SLIM and AKAP stencils,
and so their fill-in, match CM's and PN's; see
Figure~\ref{fig:sparsity_pattern}. SLIM and AKAP thus require the
same per-iteration compute cost and storage for linear solutions
as PN and CM without the same degree of convergence
benefit~\cite{Shtengel:2017:GOV}.

In summary, for smaller systems Newton-type methods have been, till
now, our likely best choice for geometry optimization, while as we
scale we have inevitably needed to move to first-order methods to
remain tractable, while accepting reduced convergence rates and
even the possibility of nonconvergence altogether.  We develop a
new quasi-Newton method, BCQN, that locally blends gradient information
with the matrix Laplacian at each iterate to regain improved and
robust convergence with efficient per-iterate storage and computation
across scales while avoiding the current pitfalls of L-BFGS methods.

\subsection{Line search}
\label{sec:rel_line_search}

Once we have applied the effort to compute a search direction we would like to maximize its effectiveness by taking as large a step along it as possible. Because the energies we treat are nonlinear, too large a step size will actually make things worse by accidentally increasing energy. A wide range of line-search methods are thus employed that search along the step direction for \emph{sufficient} decrease~\cite{Nocedal:2006:Book}. However, when we seek to minimize nonconvex energies on meshes the situation is even tougher. Most (although not all) popular and important nonlinear energies, both in geometry processing and physics, are composed by the sum of rational fractions of singular values of the deformation gradient $W(F) = W(\sigma) = f(\sigma)/g(\sigma)$ where the denominator $g(\sigma) \rightarrow 0$ as $\sigma_i  \rightarrow 0, \forall i \in [1,d]$. Notice that these $1/g(\sigma)$ barrier functions block element inversion. 
Irrespective of their source, these blocking nonconvex energies rapidly increase energy along any search direction that would collapse elements. To prevent this (and likewise the possibility of getting stuck in an inverted state) search directions are capped to prevent inversion of every element in the mesh. This is codified by Smith and Schaeffer's~\shortcite{Smith:2015:BPW} line-search filter, applied before traditional line search, that computes the maximal step size that guarantees no inversions anywhere. 

Unfortunately, this has some serious consequences for progress. Notice that if even a single element is close to inversion this can amputate the full descent step so much that almost no progress can be made at all; see Figure~\ref{fig:blocked_line_search}. This in many senses seems unfair as we should expect to be able to make progress in other regions where elements may be both far from inversion and yet also far from optimality. 
To address these barrier issues we develop an efficient barrier-aware filter that allows us to avoid blocking contributions from individual elements close to collapse while still taking large steps elsewhere in the mesh, see Figure~\ref{fig:blocked_line_search}, top.

\subsection{Termination}
\label{sec:termination_woes}

Naturally we want to take as few iterates as possible while being sure
that when we stop, we have arrived at an accurate solution according to
some easily specified tolerance. The gold-standard in optimization
is to iterate until the gradient is small $\| \nabla E \| < \epsilon$, for
a specified tolerance $\epsilon>0$. This is robust as $\nabla E$ is zero only at stationary
points, and with a bound on Hessian conditioning near the solution can even provide
an estimate on the distance of $x$ to the solution.

\begin{wrapfigure}{r}{0.5\linewidth}
  \begin{center}
    \includegraphics[width=1\linewidth]{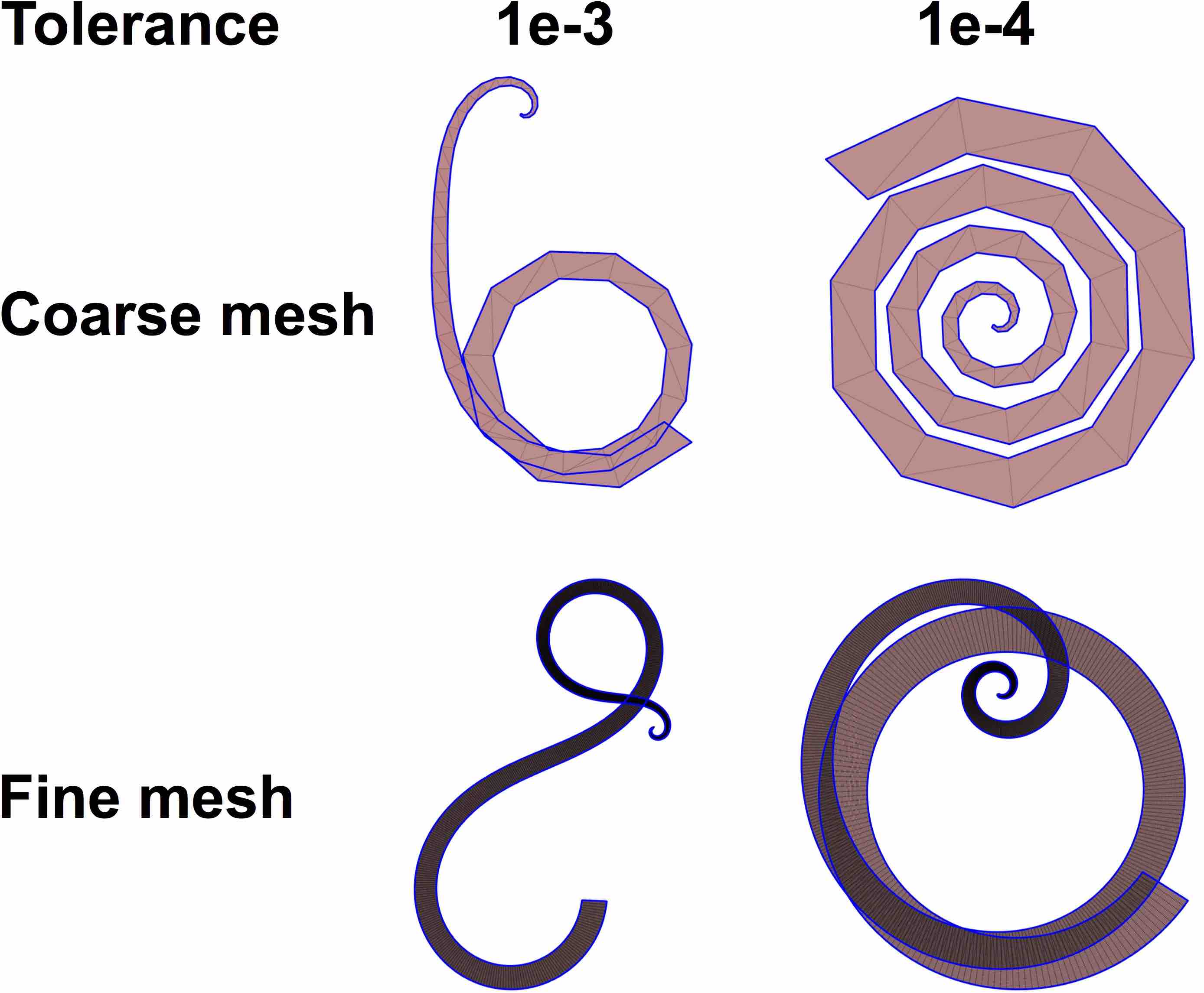}
    \caption{Standard termination measures, e.g.\ the vertex-scaled
    gradient norm above, are inconsistent across mesh, energy and scale changes.}
    \label{fig:term_compare_1}
  \end{center}
\end{wrapfigure}
However, an appropriate
value of $\epsilon$ for a given application is highly depend on the mesh, its
dimensions, degree of refinement, energy, etc.
A common engineering rule of thumb to deal with refinement consistency is to instead divide the
$L2$-norm of $\nabla E$ by the number of mesh vertices.  However, as we see in the
inset figure, this normalization does not help significantly, for
example here across changes in mesh resolution for the 2D swirl test;
see Section~\ref{sec:term_results} for more experiments.

\begin{figure}[h]
\centering
\vspace{16pt}
\includegraphics[width=1\linewidth]{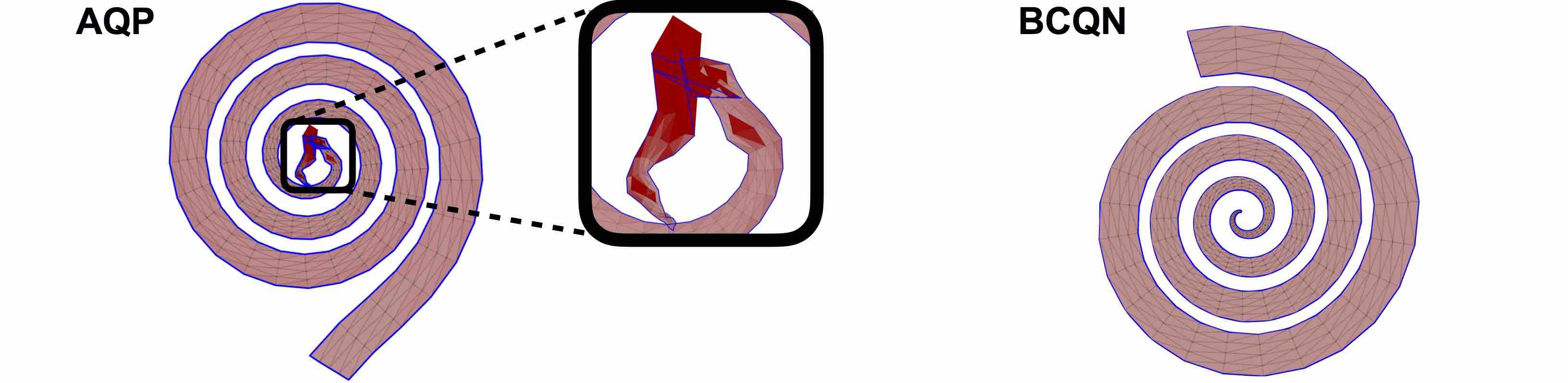}
\caption{In the 2D swirl example, BCQN with our reliable termination criterion
(\textbf{right}) only stops once it has actually reached a satsifactory solution.
The reference AQP implementation (\textbf{left}) erroneously declares success
early on when it finds two iterates have barely changed, but this is due only
to hitting a difficult configuration where AQP struggles to make progress.}
\label{fig:aqp_stop}
\end{figure}

To avoid problem dependence, recent geometry optimization
codes generally either take a fixed (small) number of
iterations~\cite{Rabinovich:2016:SLI} or iterate until an
absolute or relative error in energy $\|E_{i+1} -E_i\|$ and/or
position $\| x_{i+1} - x_i \|$ are small~\cite{Shtengel:2017:GOV,Kovalsky:2016:AQP}.
However, experiments underscore there is not yet any method
which always converges satisfactorily in the same fixed number of
iterations, no matter varying boundary conditions, shape difficulty, mesh
resolution, and choice of energy. Measuring the change in
energy or position, absolutely or in relative terms, unfortunately
cannot distinguish between an algorithm converging and simply
stagnating in its progress far from the solution; again, there is not
yet any method which can provably guarantee any degree of progress at every
iterate before true convergence. 
Figure~\ref{fig:aqp_stop} illustrates, on the swirl example, how the
reference AQP implementation declares convergence well before it reaches
a satsifactory solution, when early on it hits a difficult configuration
where it makes little local progress.

To provide reassuring termination criteria in practice and to enable
fair comparisons of current and future geometry
optimization problems we develop a gradient-based stopping criterion
which remains consistent for optimization problems even as we vary
scale, mesh resolution and energy type. This allows us, and future users,
to set a default convergence tolerance in our solver once and leave it
unchanged, independent of scale, mesh and energy. This likewise
enables us to compare algorithms without the false positives
given by non-converged algorithms that have halted due to lack of progress.

\section{Blended Quasi-Newton}
\label{sec:blend}

In this section we construct a new quadratic energy proxy which
effectively blends the Sobolev gradient with L-BFGS-style updates
to capture curvature information, avoiding the troubles previous
quasi-Newton methods have encountered in geometry optimization.
Apart from the aforementioned issue of a dense proxy incorrectly
coupling distant vertices in L-BFGS and SL-BFGS, we also find that
the gradients for non-convex energies with barriers can have highly disparate
scales, causing further trouble for L-BFGS. The much smoother
Sobolev gradient diffuses large entries from highly distorted
elements to the neighborhood, giving a much better scaling.
The Laplacian also provides essentially the correct structure for
the proxy, only directly coupling neighboring elements in the mesh,
and is well-behaved initially when far from the solution, thus we
seek to stay close to the Sobolev gradient, as much as possible, while
still capturing valuable curvature information from gradient history.

The standard (L-)BFGS approach exploits the secant approximation
from the difference in successive gradients, 
$y_i = \nabla E(x_{i+1}) - \nabla E(x_{i})$ compared to the
difference in positions $s_i = x_{i+1}-x_i$,
\begin{equation}
\label{eq:proxy_1}
\begin{aligned}
 \nabla^2 E(x_{i+1}) s_i & \simeq  y_i \\
\Rightarrow \quad \nabla^2 E(x_{i+1})^{-1} y_i & \simeq s_i,
\end{aligned}
\end{equation}
updating the current inverse proxy matrix $D_i$ (approximating
$\nabla^2 E^{-1}$ in some sense) so that $D_{i+1}y_i = s_i$.
The BFGS quasi-Newton update is generically
\begin{equation}
\label{eq:BFGS_update}
\mathrm{QN}_i(z, D) = V_i(z)^T D V_i(z) + \frac{s_i s_i^T}{s_i^Tz},  \> \> V_i(z) = I - \tfrac{z s_i^T}{s_i^Tz}.
\end{equation}
We can understand this as using a projection matrix $V_i$ to annihilate
the old $D$'s action on $z$, then adding a positive semi-definite
symmetric rank-one matrix to
enforce $\mathrm{QN}_i(z,D)z = s_i$. Classic BFGS uses
$D_{i+1} = \mathrm{QN}_i(y_i, D_i)$, whereas L-BFGS uses
\begin{equation}
    D_{i+1} = \mathrm{QN}_i(y_i, \tilde{D}_i),
\end{equation}
where $\tilde{D}_i$ has the oldest $\mathrm{QN}$ update removed,
and crucially represents each $D$ as a product of linear operators,
rather than an explicit full matrix. Only the last $m$ $\{s,y\}$ vector pairs (we
use $m=5$) along with the initial $D_1$ (we use the inverse Laplacian,
storing only its Cholesky factor) are stored; application of $D$ is
then just a few vector dot-products and updates along with backsolves for
the Laplacian.

\begin{figure}[h!]
\centering
\includegraphics[width=0.9\linewidth]{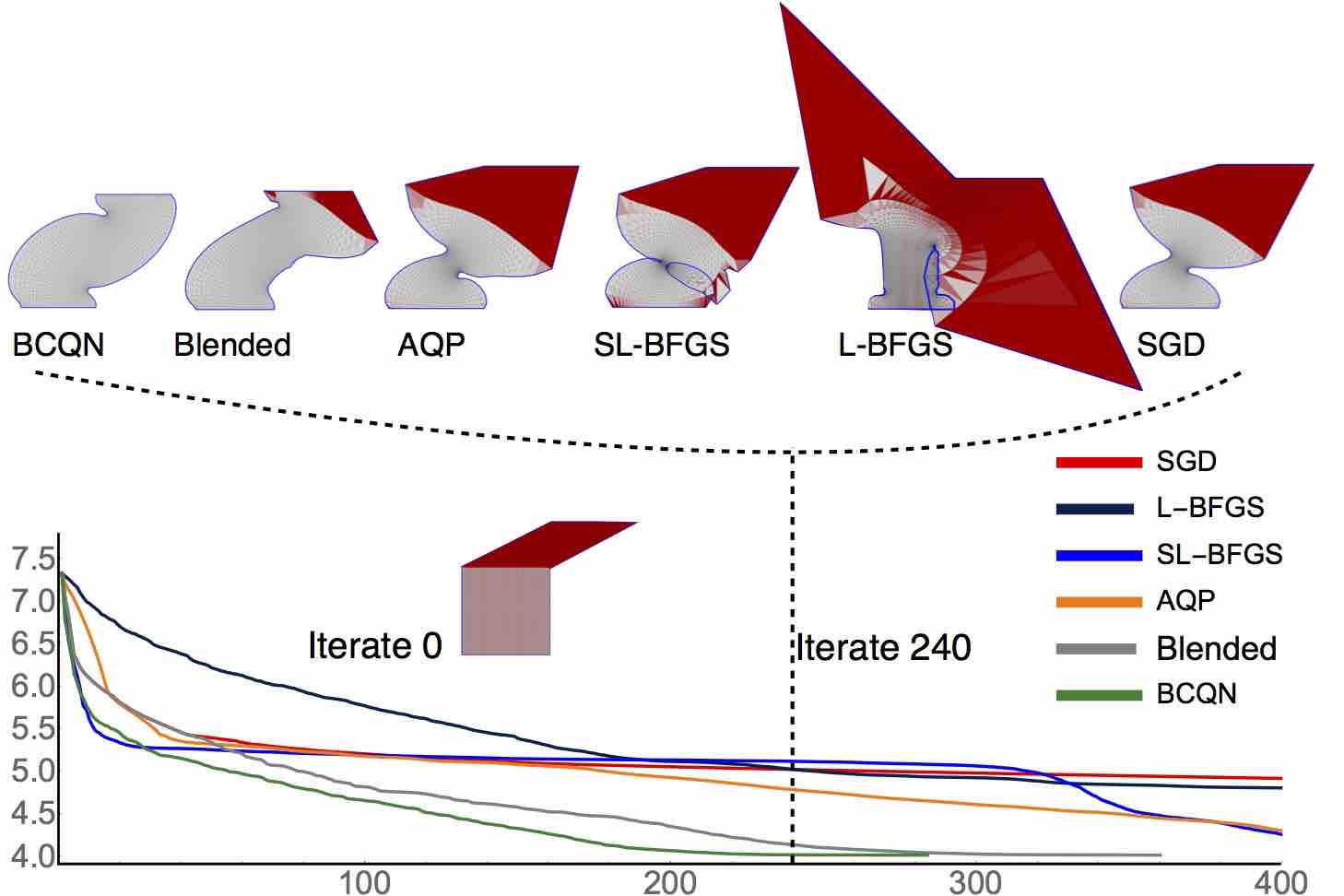}
\caption{A 2D shearing deformation stress
test with MIPS energy, comparing methods by plotting iteration vs.\ energy. Both L-BFGS 
as well as inverse Laplacian initialized (SL-BFGS) have slow convergence as previously
reported -- especially when compared to SGD and AQP which use just
the Laplacian. 
At iteration 240 the visualized deformations show both L-BFGS-based methods suffering from
swelling due to inaccurate coupling of distant elements.
Applying our blending model alone (Blended) is highly
effective, while our full BCQN method gives the best results overall.}
\label{fig:quadratic_compare}
\end{figure}

\subsection{Greedy Laplacian Blending}

Experiments show that far from the solution, the Laplacian is often a much
more effective proxy than the L-BFGS secant version: see AQP/SGD vs.\ L-BFGS in Figure \ref{fig:quadratic_compare}.
In particular, the difference in energies $y$ may introduce spurious coupling or
have badly scaled entries near distorted triangles. In this
case if the energy were based on the Laplacian itself (the 
\emph{Dirichlet} energy), the difference in gradients would be the better
behaved $Ls$. This motivates trying the update with $Ls$ instead of $y$,
\begin{equation}
\label{eq:qn_L}
D_{i+1} = \mathrm{QN}_i(L s_i, \tilde{D}_i),
\end{equation}
which will keep us consistent with Sobolev preconditioning, which is very effective
in initial iterations. However, to achieve the
superlinear convergence L-BFGS offers, near the solution we will want
to come closer to satsifying the secant equation, switching to using $y$ instead.

We can thus imagine a blending strategy, which uses
\begin{equation}
z_i=(1-\beta_i)y_i + \beta_i Ls_i
\end{equation}
in $\mathrm{QN}(z_i, \tilde{D}_i)$, with blending parameter $\beta_i \in [0,1]$.
A greedy strategy might choose $\beta_i$ to scale $Ls_i$ to be as close to $y_i$ as possible,
\begin{equation}
\label{eq:BQCN_proj}
\beta_i = \argmin_{\beta\in[0,1]} \| y_i - \beta L s_i \|^2,
\end{equation}
in other words using the projection of $y_i$ onto $Ls_i$. This comes as close as possible to
satsifying the secant equation with $Ls_i$, then makes up the rest with $y_i$. 
Solving (\ref{eq:BQCN_proj}) gives 
\begin{align}
\label{eq:BCQN_proj2}
\beta_i = \mathrm{proj}_{[0,1]} \left( \frac{{y_i}^T L s_i}{ \|L s_i\|^2 } \right).
\end{align}
Observe that when $Ls$ is roughly aligned with the gradient jump $y$ , but $y$ is as large or larger, $\beta$ grows and Laplacian smoothing increases --- as we might
hope for initially when far from the solution, where the Sobolev gradient is most effective.
When the energy Hessian diverges strongly from from the Laplacian
approximation, perhaps when the cross-terms between coordinates missing from the scalar Laplacian
are important, then $\beta$ will decrease, so that contributions from $y_i$ again grow.
Finally, as the gradient magnitudes decreases close to the solution,
$\beta$ will similarly decay, ideally regaining the
superlinear convergence of L-BFGS near local minima.

\subsection{Blended Quasi-Newton}

With the blending projection (\ref{eq:BCQN_proj2}) in place we experimented with a range of rescalings in hopes of 
further improving efficiency and robustness. After extensive testing we have so far found the following scaling
to offer the best performance:
\begin{align}
\begin{split}
\beta_i =  \mathrm{proj}_{[0,1]} \Big(\frac{ \mathrm{normest}(L) {y_i}^T L s_i}{A(V,T)} \Big), \\
\text{with} \> \> A(V,T) = \Big(\sum_{t \in T} a_t \Big)^{\frac{2 (d - 1)}{d}}.
\end{split}
\end{align}
Here $\mathrm{normest}(L)$ is an efficient estimate of the matrix 2-norm using power iteration,
and $A(V,T)$ is a \emph{constant} normalizing term with appropriate dimensions and so no longer
has the same potential concern for sensitivity in the denominator when $Ls$ is small but $s$ isn't. Both terms are computed just
once before iterations begin and reused throughout. 

As mentioned, we initialize the inverse proxy with $D_1=L^{-1}$,
thus starting with Laplacian preconditioning. With line search
satisfying Wolfe confitions our proxy remains SPD across all
steps~\cite{Nocedal:2006:Book}. Each step jointly updates $D_i$
using the standard two-loop recursion and finds the next descent
direction $s_i = -D_i \nabla E(x_i)$.
Figure~\ref{fig:quadratic_compare} illustrates the gains possible from
blended quasi-Newton compared to both standard L-BFGS and Sobolev gradient algorithms, while then applying our barrier-aware filter, derived in our next section gives best results with our
full BCQN algorithm.

\section{Barrier-Aware Line Search Filtering}
\label{sec:constraints}

As mentioned in Section~\ref{sec:rel_line_search} and shown
in Figure~\ref{fig:blocked_line_search}, the barrier factor $1/g(\sigma)$
in nonconvex energies typically dominates step size in line search.
Even a single element that is brought close to collapse by the
descent direction, $p_i$, can restrict the line search step size
severely.  The computed step size $\alpha_i$ then scales $p_i$
\emph{globally} so that all elements, not just those that are going
to collapse along $p_i$, are prevented from making progress. To avoid
this, a natural strategy suggests itself: when the descent direction would cause
elements to degenerate towards collapse along the full step,
rather than simply truncating line search as in Smith and
Schaefer\ \shortcite{Smith:2015:BPW}, we filter collapsing contributions
from the search direction prior to line search.
We call this strategy \emph{barrier-aware line search filtering}.

\subsection{Curing line search}

\begin{figure}[t!]
\centering
\includegraphics[width=1\linewidth]{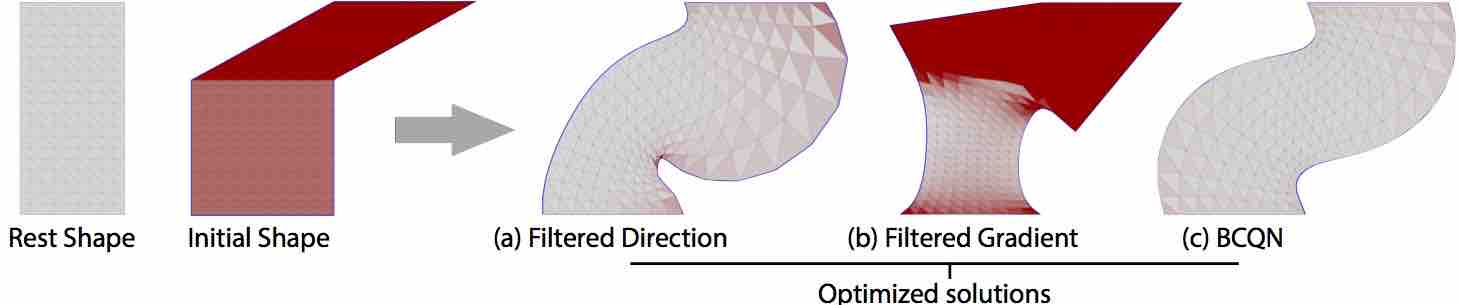}
\caption{
\bfi{Direct filtering does not work.} Zeroing out inverting components of descent directions
or gradients makes the search direction inconsistent with the objective and so prevents convergence,
leading to termination at poor solutions (a) and (b).
{\bf Left:} we initialize a 2D shear deformation,
constraining the top of a bar to slide rightwards.
{\bf Middle:} direct filtering of the descent direction (a) and
the gradient (b) allow large descent steps forward unblocked
from the contributions of close-to-collapsed elements. However, this
results in termination at shapes that that do not satisfy optimality of the original minimization.
{\bf Right:} compare to an optimal solution for this problem (c)
obtained with BCQN.
 }
\label{fig:filter_fail}
\end{figure}

Figure~\ref{fig:filter_fail} illustrates how the simplest possible filters,
zeroing out contributions from nearly-inverted elements
in either the search direction (\ref{fig:filter_fail}a)
or the gradient before Laplacian smoothing (\ref{fig:filter_fail}b)
fail. We must be able to make progress in nearly-inverted elements
when the search direction can help, or there is no hope for reaching
the actual solution; simple zeroing fails to converge, which is
no surprise as it in essence is arbitrarily manipulating the
target energy, changing the problem being solved.
We instead want to \emph{augment} the original optimization problem
in a way which doesn't change the solution, but gives us a tool to
safely deal with problem elements so the search direction $p_i$ doesn't
cause them to invert, ideally with a small fixed cost per iteration.

\subsection{One-Sided Barriers in Geometry Optimization}

Element $t \in T$ is inverted at positions $x$ precisely when the orientation function
$a_t(x) = \det(F_t(x))$ is negative. Concatenating over $T$, the global vector-valued function for element
orientations is then
\begin{equation}
a(\cdot) = \big(a_1(\cdot), ..., a_m(\cdot) \big)^T.
\end{equation}
As long as $a(x) > 0$, no element is collapsed or inverted, and the energy remains finite.
Note, however, many energies are also finite for inverted elements $a_t(x)<0$, only blowing up
at collapse $a_t(x)=0$, so technically there may exist local minima where $\nabla E(x^*)=0$
yet some elements are inverted. Generally, practitioners wish to rule these potential solutions
out however, with two implicit but so far informal assumptions of locality: 
the initial guess is not inverted, $a(x_1)>0$, and that the solver follows a path
which never jumps through the barrier to inversion. 

We formalize these requirements in the optimization as
\begin{equation}
\label{eq:hard_constr_E}
\min_x \{E(x) \ : \  a(x) \geq 0 \}.
\end{equation}
Adding the constraint $a(x) \geq  0$ now explicitly restricts our
optimization to noninverting deformations but otherwise leaves the
desired solution unchanged. (See Supplement, Section 1, for proof.)

\subsection{Iterating Away from Collapse}

With problem statement (\ref{eq:hard_constr_E}) in place, we now exploit it in curing the search direction from
collapsing elements. At each iterate $i$, form the projection  
\begin{align}
\label{eq:p_project}
 \min_p \left\{ \| p + D_i \nabla E(x_i) \|_2^2 \> : \> a(x_i) + \nabla a(x_i)^T p \geq 0 \right\} 
\end{align}
of the predicted descent direction $\tilde{p}_i = -D_i \nabla E(x_i)$ onto 
the subset satisfying a linearization of the no-collapse condition.
Satisfying (\ref{eq:p_project}) exactly would ensure that projected
directions would not locally generate collapse and likewise preserve
symmetry~\cite{SKVTG2012}. However, its exact solution is neither
necessary nor efficient. Instead, we construct an approximate
solution to (\ref{eq:p_project}) as a filter that \emph{helps}
avoid collapse, preserves symmetry, and guarantees a low cost for
computation for all descent steps. 

Strict convexity of the projection guarantees that a minimizer $p^*$ of (\ref{eq:p_project}) is given by the
KKT\footnote{Here and in the following $\lambda = (\lambda_1, ..  ,\lambda_m)^T \in R^m$
is a Lagrange multiplier vector and $\vc x \perp \vc y$ is the \emph{complementarity condition}
$y_t z_t = 0,\ \forall t$.} conditions~\cite{Bertsekas:2016:NOP}
\begin{align}
\label{eq:kkt_prog1}
p^*+ D_i \nabla E(x_i) - \nabla a(x_i) \lambda^* = 0, \\
\label{eq:kkt_prog2}
0 \leq \lambda^* \perp a(x_i) + \nabla a(x_i)^T p^* \geq 0.
\end{align}
We simplify with $C_i = \nabla a(x_i)$, $M_i = \nabla a(x_i)^T \nabla a(x_i)$,
and $b_i = a(x_i)$, then form the Schur complement of the above to arrive at an equivalent
Linear Complementarity Problem (LCP)~\cite{Cottle:2009}
\begin{align}
\label{eq:LCP_proj}
\begin{split}
0 \leq \lambda^* \perp M_i \lambda^* + C_i^T  p_i + b_i \geq 0,
\end{split}
\end{align}
and then a damped Jacobi splitting
with $M_i = \omega^{-1}  T_i +  (M_i - \omega^{-1} T_i)$,
diagonal $T_i = \mathrm{diag}(M_i)$ and damping parameter
$\omega \in (0,1)$. This gives us an iterated LCP ranging over
iteration superscripts $j$,
\begin{align}
\label{eq:LCP_proj_split}
\begin{split}
0 \leq \lambda^{j+1} \perp \omega^{-1} T_i \lambda^{j+1} + M_i \lambda^j - \omega^{-1} T_i \lambda^j + C_i^T  p_i  + b_i \geq 0.
\end{split}
\end{align}

\begin{figure}[t!]
\centering
\includegraphics[width=0.9\linewidth]{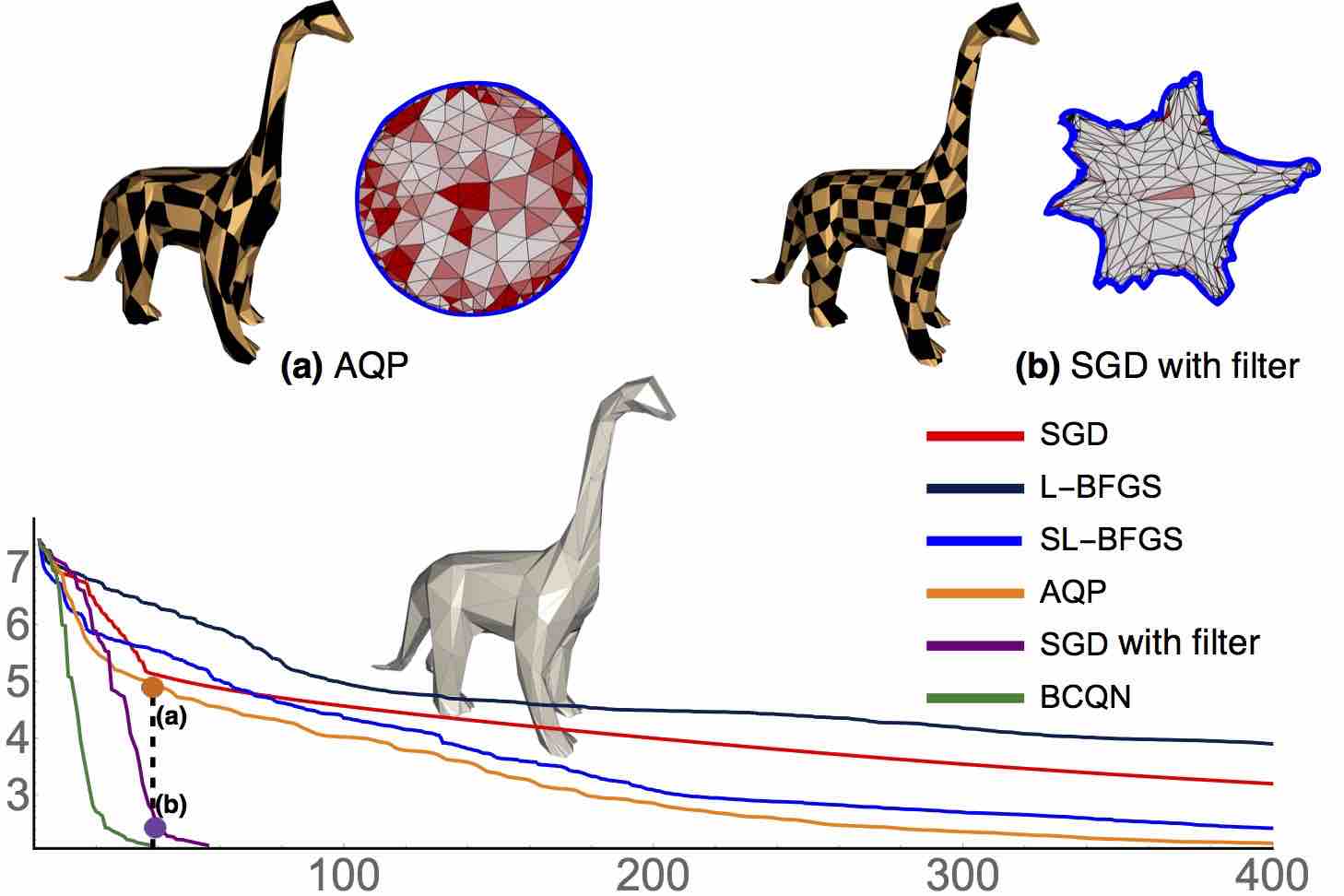}
\caption{
\bfi{Line search filtering.} {\bf Bottom:} We optimize a
uv-parameterization with the MIPS energy to consider line search
filtering behavior, plotting energy (y-axis) against iteration
counts for a range of methods. Just adding our barrier-aware line
search filtering alone to SGD improves its convergence by
well over an order of magnitude, and almost an order of magnitude
over AQP as well as plain L-BFGS and SL-BFGS. BCQN with blending
and line search filtering improves convergence even further.
{\bf Top:} a comparison of the embeddings and texture-maps
for AQP and SGD with the filter at the $40^\textrm{th}$
iterate.
}
\label{fig:combined_method}
\end{figure}

\subsection{Line Search Filtering}

Each iteration of the splitting (\ref{eq:LCP_proj_split}) simplifies
to the damped projected Jacobi (DPJ) update\footnote{We use the convention $[\cdot]^+ = \max[0, \cdot]$.}
\begin{align}
\label{eq:DPJ}
\lambda^{j+1} \leftarrow \left[\lambda^j - \omega T^{-1}\big(C_i^T (C_i \lambda^j) + c_i\big)\right]^+,
\end{align}
with constant $c_i = C_i^T  p_i + b_i$. Here each of the $m$ entries in $\lambda^{j+1}$
can be updated in parallel (unlike with Gauss-Seidel iteration).
As $M_i$ is PSD this iteration process converges to
(\ref{eq:LCP_proj})~\cite{Cottle:2009} and so to (\ref{eq:p_project}).
We do not seek a tight solution, however, as we just want to be sure the worst blocks
to line search are filtered away. Therefore we initialize with $\lambda^0=0$ to avoid
unnecessary perturbation, use a coarse termination tolerance
for DPJ (see below), and never use more than a maximum of 20 DPJ iterations.

At each DPJ iteration $j$ we check for termination with an LCP
specialized measure, the Fischer-Burmeister
function~\cite{Fischer:1992:ASN}
$\mathrm{FB}(\lambda^j, M_i  \lambda^j  +  c_i)$ evaluated as
\begin{align}
\label{eq:FB}
\mathrm{FB}(a,b) = \sqrt{\sum_{k \in [1,m]}  \left(a_k + b_k - \sqrt{a_k^2 + b_k^2} \right)^2}.
\end{align} 
As we initialize with $\lambda^0 = 0$, when $p_i$ is non-collapsing
$\mathrm{FB} = 0$, and thus no line search filtering iterations
will be applied. Likewise, we stop iterations whenever the $\mathrm{FB}$
measure is roughly satisfied by either a relative error of $<10^{-3}$
or an absolute error $<10^{-6}$.

Filtering thus applies a fixed maximum upper limit on computation
and performs no iterations when not necessary. Upon termination of
DPJ iterations, plugging our final $\lambda$ into (\ref{eq:kkt_prog1})
we obtain our update to form the line search filtered descent
direction
\begin{align}
p^\ell_i = p_i  + C_i \lambda.
\end{align}
As Figure~\ref{fig:blocked_line_search} shows, despite the rough
nature of the filter, it can make a dramatic difference in line search.

\section{Termination Criteria}
\label{sec:term}

Every iterative method for minimizing an objective function $E(x)$
must incorporate stopping criteria: when should an approximate
solution be considered good enough to stop and claim success?
Clearly, in the usual case where the actual minimum value of $E(x)$
is unknown, basing the test on the current value of $E(x_i)$ is
futile. As noted in Section~\ref{sec:termination_woes},
stopping when successive iterates are closer than some tolerance is
vulnerable to false positives (halting far from a solution), as is
using a fixed number of iterations. Although monitoring $\|\nabla E\|$
is robust, each individual problem may need a different tolerance to
define a satisfactory solution even when normalized by number of vertices:
see Figures\ \ref{fig:term_compare_1} and\ \ref{fig:term_compare_2}.
We thus propose a new way to derive and construct an appropriate,
roughly problem-independent, relative scale for a gradient-based
measure for a stopping criterion.

\subsection{Characteristic Gradient Norm}

All energies we consider are summations of per-element energy densities $W(\cdot)$ computed
from the deformation gradient $F_t(x)$ and weights $a_t$, in each element $t$, as per equation (\ref{eq:obj}). 
To simplify the following we can then evaluate energy densities on
the vectorized deformation gradient as $W\big(vec(F_t)\big) =  W(G_t
x)$, where $G_t$ is the linear gradient operator for element $t$.
The full energy gradient is then
\begin{equation}
    \nabla E(x) = \sum_{t \in T} a_t  G_t^T \nabla W(G_t x).
\end{equation}
We wish to generate a ``characteristic'' value we can compare this gradient to meaningfully, with the same
dimensions; we will do this with each component of the above summation separately.

First observe that the deformation gradient, $F_t$, the argument to $W$, is dimensionless and therefore
$\nabla W$ has the same dimensions as $W$, and even as the element Hessian $\nabla^2 W$. For the simplest
quadratic energy densities, this Hessian has the attractive property of being constant; we thus choose
to use the 2-norm of the Hessian, evaluated about the deformation gradient at rest ($F_t=I$), to get
a representative value for $\nabla W$:
\begin{equation}
    \langle W \rangle = \|\nabla^2 W(I)\|_2.
\end{equation}

Second, note that the $i^\textrm{th}$ part of $G_t$ for a triangle
(respectively tetrahedra) $t$ containing vertex $i$ will attain its
maximum value for fields which are constant along the opposing edge
(triangle) and that value will be the reciprocal of the altitude.
Up to a factor of $2$ ($3$), this is the length (area) of the
opposing edge (triangle) divided by the rest area (volume), of the
element, i.e.\ $a_t$. Summing over all incident elements, weighted
by $a_t$, we arrive at a characteristic value for vertex $i$ of
$\ell_i$ equalling the perimeter (surface) area of the one-ring of
vertex $i$. We compute this value for all vertices, giving us the
vector $\ell(V,T) = (\ell_1, ..., \ell_n)^T \in \R^n$, with one
scalar entry per vertex.

The product of our energy and mesh values together form the characteristic value for the norm of the gradient
\begin{equation}
    \langle W \rangle \| \ell(V,T) \|,
\end{equation}
where we take the same vector norm as that with which we evaluate $\|\nabla E(x)\|$; we use the 2-norm in all our
experiments. For all methods we stop iterating when
\begin{equation}
    \|\nabla E(x)\| \leq \epsilon \langle W \rangle \| \ell(V,T) \|,
\end{equation}
given a dimensionless tolerance $\epsilon$ from the user, which is
now essentially mesh- and energy-independent. See Figures\
\ref{fig:term_compare_1}, \ref{fig:ours_yaron} and\ \ref{fig:term_compare_2} as
well as our experimental analysis in Section\ \ref{sec:results} for evaluation.

\section{The BCQN Algorthim}
\label{sec:alg}

\begin{algorithm}[h!]
\label{alg:BCQN}
\caption{Blended Cured Quasi-Newton (BCQN)}

\textbf{Given:} $x_1$, $E$, $\epsilon$  \hspace{10pt} 

\textbf{Initialize and Precompute:}

\hspace{30pt} $s = \epsilon \langle W \rangle \| \ell(V,T) \| $  \hspace{10pt} // Characteristic termination value (\S\ref{sec:term})

\hspace{30pt} $L, \> \> D \leftarrow L^{-1}$ \hspace{10pt} // Initialize blend model  (\S\ref{sec:blend})

\hspace{30pt} $g_1 = \nabla E(x_1), \> \> i = 1$ 
 
$\textbf{while}$ $ \|g_i\| > s$ $\textbf{do}$\hspace{10pt}// Termination criteria  (\S\ref{sec:term})\\

\hspace{10pt} $p \leftarrow -D g_i$ \hspace{10pt}//  Precondition gradient (\S\ref{sec:blend})\\

\hspace{10pt} // Assemble for DPJ iterations (\S\ref{sec:constraints}):

\hspace{20pt} $C \leftarrow \nabla a(x_i)$ 

\hspace{20pt} $M \leftarrow C^T C, \> \> c \leftarrow C^T  p + a(x_i)$

\hspace{20pt} $E \leftarrow \mathrm{diag}(M)^{-1}, \> \> \lambda \leftarrow 0$ 

\hspace{10pt} $\mathit{fb} \leftarrow \textrm{FB}(\lambda, M \lambda  +  c)$ \hspace{10pt}// LCP residual (Equation (\ref{eq:FB}) in \S\ref{sec:constraints})

\hspace{10pt} \textbf{for} $j = 1$ \textbf{to} \text{20} \hspace{10pt}// Line-search preconditioning  (\S\ref{sec:constraints})

\hspace{20pt}\textbf{if} $\mathit{fb} < 10^{-6}$ \textbf{then} \hspace{3pt} \textbf{break} \hspace{5pt} \textbf{end if}

\hspace{20pt}$\mathit{fb} \leftarrow \mathit{fb}_\mathrm{next}$ 

\hspace{20pt} $\lambda \leftarrow [\lambda - \tfrac{1}{2} E \big(C^T (C \lambda) + c\big)]^+$ // Parallel project  (\S\ref{sec:constraints})

\hspace{20pt}$\mathit{fb}_\mathrm{next} \leftarrow \textrm{FB}(\lambda, M  \lambda  +  c)$  

\hspace{20pt}\textbf{if}  $|\mathit{fb} - \mathit{fb}_\mathrm{next}| / \mathit{fb} < 10^{-3}$ \textbf{then} \hspace{3pt} \textbf{break} \hspace{5pt} \textbf{end if}

\hspace{10pt}\textbf{end for}

\hspace{10pt} $p^\ell \leftarrow p + C \lambda$ \hspace{10pt}// Line-search filtered search direction  (\S\ref{sec:constraints})

\hspace{10pt}$\alpha \leftarrow \text{LineSearch}(x_i, p^\ell, E)$ \hspace{10pt} // Line search (\S\ref{sec:blend})

\hspace{10pt}$x_{i+1} = x_i + \alpha p^\ell$  \hspace{10pt} // Descent step (\S\ref{sec:blend})

\hspace{10pt} $g_{i+1} = \nabla E(x_{i+1})$

\hspace{10pt}$D \leftarrow \text{Blend}(D, L, x_{i+1}, x_i, g_{i+1}, g_i)$\hspace{10pt}// BCQN blending update (\S\ref{sec:blend})

\hspace{10pt}$i \leftarrow i+1$

$\textbf{end while}$

\end{algorithm}


\begin{figure}
\centering
\includegraphics[width=1\linewidth]{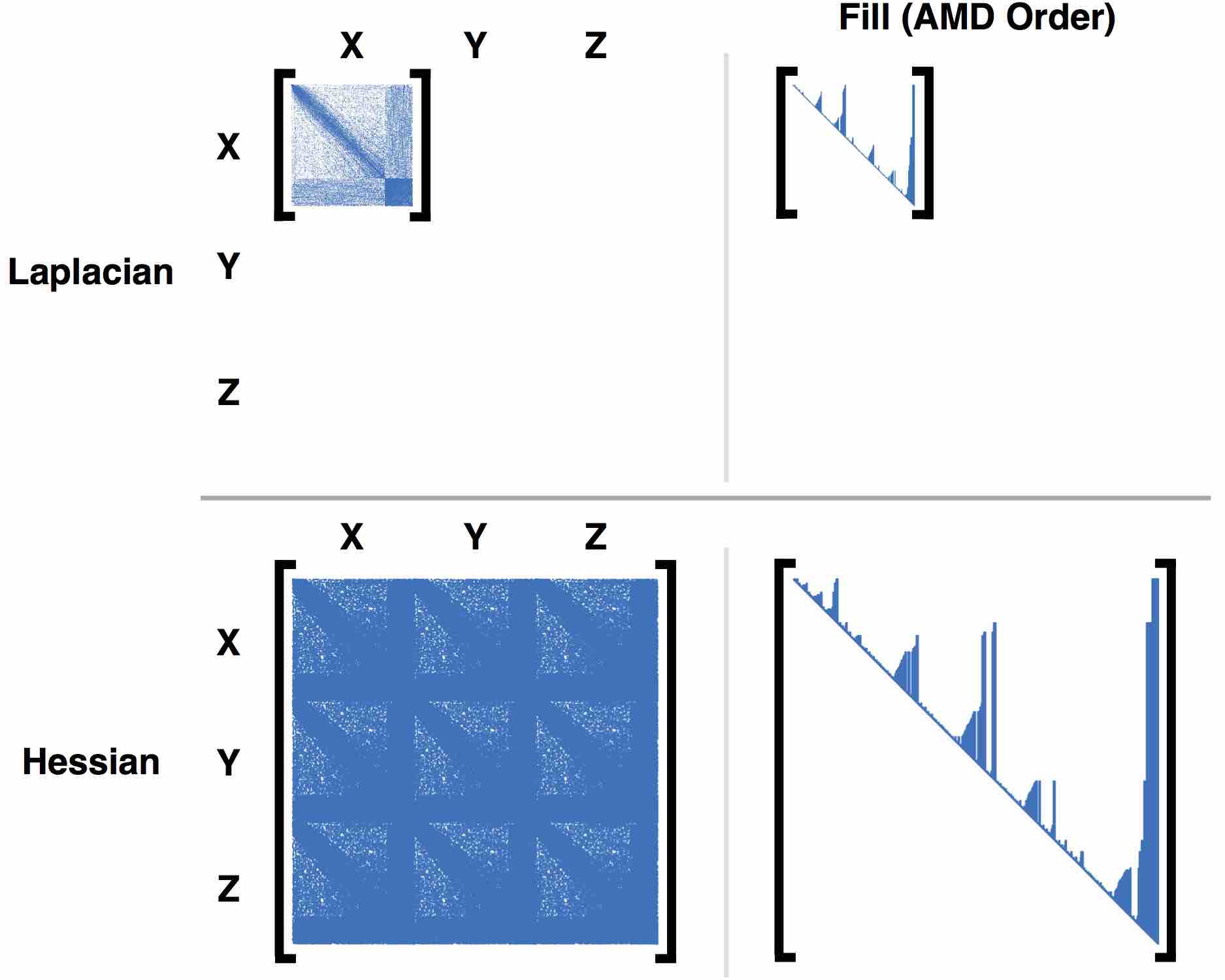}
\caption{\bfi{Sparsity Differences in Proxies.} {\bf Left:} The scalar Laplacian (top) is smaller \emph{and} sparser than the
Hessian and its approximations (bottom) used in CM, PN, SLIM and AKAP. {\bf Right:} This results in a much cheaper factorization and solve for the Laplacian; it
is applied in both BCQN and AQP independently to each coordinate.}
\label{fig:sparsity_pattern}
\end{figure}

Algorithm 1 contains our full BCQN algorithm in pseudocode. The
dominant cost, for both memory and runtime, is the Laplacian solve
embedded in the application of $D$, which again is not stored as a
single matrix, but rather is a linear transformation involving a
few sparse triangular solves with the Laplacian's Cholesky factor
and outer-product updates with a small fixed number of L-BFGS history
vectors. Recall that we separately solve for each coordinate with
a scalar Laplacian, not using a larger vector Laplacian on all
coordinates simultaneously; this also exposes some trivial parallelism.
Apart from the Laplacian, all steps are either linear (dot-products,
vector updates, gradient evaluations, etc.) or typically sublinear
(DJP assembly and iterations, which only operate on the small number
of collapsing triangles, and again are easily parallelized).

As Lipton et al.\ proved \shortcite{lipton:1979:gnd},
the lower bounds for Cholesky factorization on a two-dimensional
mesh problem with $n$ degrees of freedom are $O(n \log n)$ space
and $O(n^{3/2})$ sequential time, and in three-dimensional problems 
where vertex separators are at least $O(n^{2/3})$, their Theorem
10 shows the lower bounds are $O(n^{4/3})$ space and $O(n^2)$
sequential time. On moderate size problems running on current
computers, the cost to transfer memory tends to dominate arithmetic,
so the space bound is more critical until very large problem sizes
are reached.

\subsection{Comparison with Other Algorithms}

The per-iterate performance profile of AQP is most similar to BCQN:
it too is dominated by a Laplacian solve. The only difference is
the extra linear and sublinear work which BCQN does for the quasi-Newton
update and the barrier-aware filtering; even on small problems, this
overhead is usually well under half the time BCQN spends, and as the
next section will show, the improved convergence properties of BCQN
render it faster.

The second-order methods we compare against, PN and CM, as well as the more approximate
proxy methods, SLIM and AKAP, all use a fuller stencil which couples coordinates.
The same asymptotics for Cholesky apply, but whereas AQP and BCQN can solve a scalar $n\times n$ Laplacian
$d$ times (once for each coordinate, independently), these other methods
must solve a single denser $nd\times nd$ matrix, with $d^2$ times more
nonzeros: see Figure~\ref{fig:sparsity_pattern}.
Moreover, the matrix changes at each iteration and must be refactored,
adding substantially to the cost: factorization is significantly slower
than backsolves.

\section{Evaluation}
\label{sec:results}

\subsection{Implementation}

We implemented a common test-harness code to enable the consistent
evaluation of the comparitive performance and convergence behavior
of SGD, PN, CM, AQP, L-BFGS and BCQN across a range of energies and
geometry optimization tasks including parameterization as well as
2D and 3D deformations, where these methods allow. For AQP this extends
the number of energies it can be tested with, while more generally
providing a consistent environment for evaluating all methods. We
hope that this code will also help support the future evaluation
and development of new methods for geometry optimization.

The main body of the test code is in MATLAB to support rapid
prototyping.  All linear system solves are performed with MATLAB's
native calls to SuiteSparse~\cite{Chen:2008:ACS} with additional
computational-heavy modules, primarily common energy, gradient and
iterative LCP evaluations, implemented in C++.  As linear solves
are the bottleneck in all methods covered here, an additional
speed-up to all methods is possible with
Pardiso~\cite{Petra:2014:AAI,Petra:2014:RTS} in place of SuiteSparse;
however, as discussed in Section~\ref{sec:pardiso} this does not change
the relative merits of the methods, and would add an additional external
dependency to the test code. For verification we also confirm that
iterations in the test-harness AQP and CM implementations match the
official AQP~\cite{Kovalsky:2016:AQP} and CM~\cite{Shtengel:2017:GOV}
codes.

All experiments were timed on a four-core Intel
3.50GHz CPU. We have parallelized the damped Jacobi LCP iterations
with Intel TBB; with more cores the overhead reported below for LCP
iterations is expected to diminish rapidly.
For all UV parameterization problems we compute 
initial locally injective embeddings via the initialization code
from Kovalsky et al.~\shortcite{Kovalsky:2016:AQP}. On rare occasions
this code fails to find a locally injective map, so we
then revert to a Tutte embedding as a failsafe using the initialization
code from Rabinovich et al.~\shortcite{Rabinovich:2016:SLI}. To
enforce Dirichlet boundary conditions, i.e.\ positional constraints,
we use a standard subspace projection~\cite{Nocedal:2006:Book}, i.e.\
removing those degrees of freedom from the problem.
When line search is employed we first find a
maximal non-inverting step size with Smith and
Schaefer's method~\shortcite{Smith:2015:BPW}, followed by standard
line search with Armijo and curvature conditions.

\subsection{Termination} 
\label{sec:term_results}

To evaluate termination criteria behavior we first instrumented two
geometry optimization stress-test examples: the \emph{Swirl}
deformation~\cite{Chen:2013:PSI} and the \emph{Hilbert curve} UV
parametrization~\cite{Smith:2015:BPW}. We run both examples to
convergence ($10^{-6}$ using our characteristic gradient) reaching
the final target shapes for each. Within these optimizations we
record the 2-norm of gradient, the vertex-normalized 2-norm of
gradient, the
relative error measure~\cite{Kovalsky:2016:AQP,Shtengel:2017:GOV}
and our characteristic gradient norm for all iterations.

\begin{figure}[t]
\vspace{3mm}
\centering
\includegraphics[width=1\linewidth]{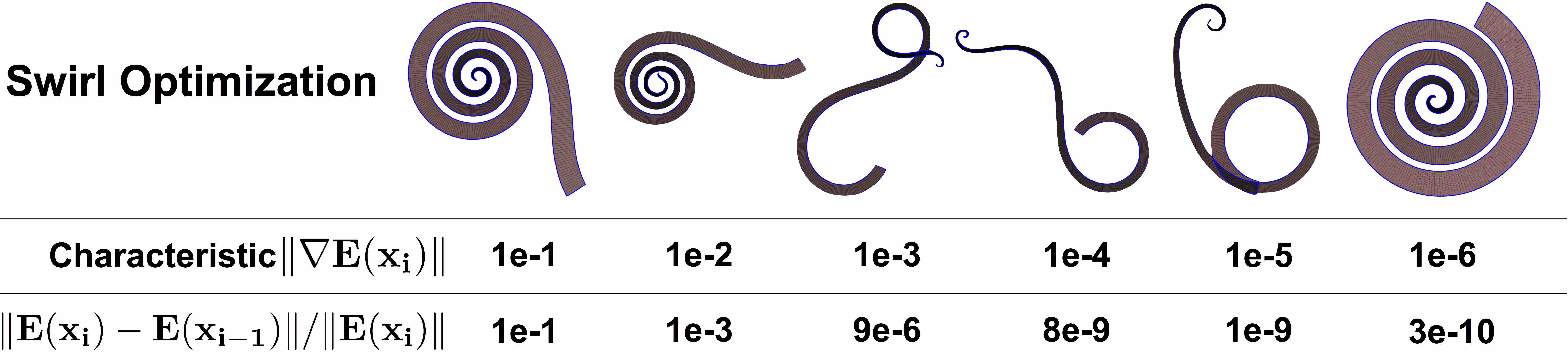}
\caption{\bfi{Termination criteria comparison.} {\bf Left to right:}
We find key points in the sequential progress of the optimized mesh
in the Swirl optimization (ISO energy) example at regular intervals of $10\times$ decrease
in our characteristic norm. We compare with the relative error measures at these same points.}
\label{fig:swirl_term_compare}
\vspace{3mm}
\end{figure}

Figure~\ref{fig:swirl_term_compare} shows the Swirl mesh
obtained during BCQN iteration at regular intervals of $10\times$
decrease in our characteristic norm. Observe
that they correspond to natural points of progress;
see our supplemental video of the entire optimization sequence
for reference. For comparison we also provide the corresponding relative
error measures, which varies much less steadily.

\begin{figure}[h!]
\vspace{3mm}
\centering
\includegraphics[width=1\linewidth]{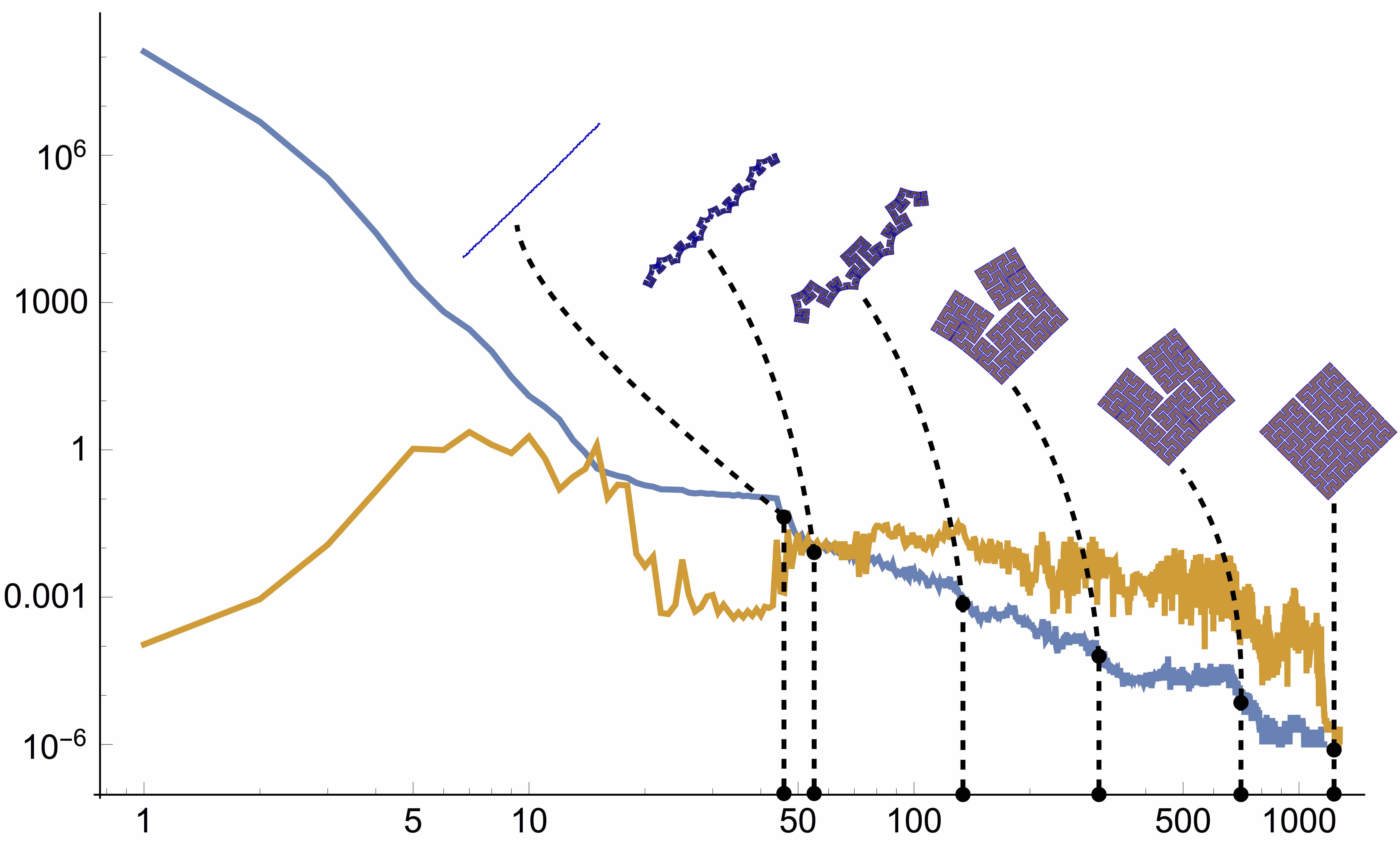}
\caption{\bfi{Measuring improvement.} 
Solving a UV parametrization of the Hilbert curve with BCQN, we plot our
characteristic gradient norm in \blue{blue} and the relative energy
error in \orange{orange} as the method proceeds, on a logarithmic scale.
Iterates are shown at decreases in the characteristic gradient norm
by factors of 10, illustrating its efficacy as a global measure of progress,
while the relative energy error measures only local changes with little
overall trend.}
\label{fig:ours_yaron}
\vspace{3mm}
\end{figure}

In Figure~\ref{fig:ours_yaron} we compare
termination criteria more closely for a UV parametrization
problem, the Hilbert curve example. We plot our characteristic gradient norm
(blue) and the relative energy error~\cite{Kovalsky:2016:AQP,Shtengel:2017:GOV}
(orange) as BCQN proceeds. Note
that the characteristic gradient norm provides consistent decrease
corresponding to improved shapes and so provides a
practical measure of improvement. The local error in energy, on the
other hand, varies greatly, making it impossible to judge how much
global progress has been made towards the optimum.

\begin{figure}[h!]
\vspace{3mm}
\centering
\includegraphics[width=1\linewidth]{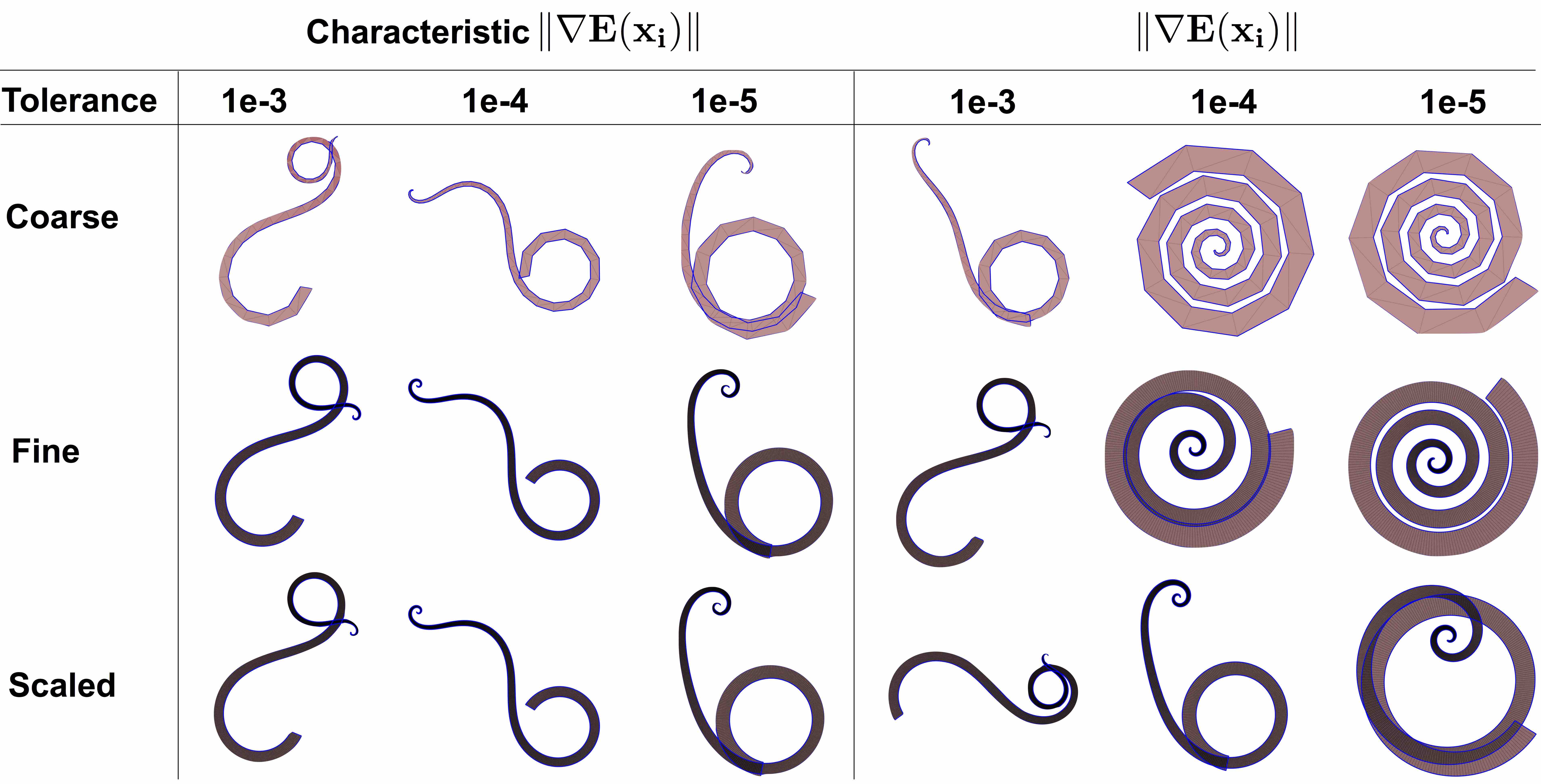}
\caption{\bfi{Termination criteria comparison across mesh refinement and scale.}
{\bf Left and right}: we show the Swirl optimization
when our characteristic norm (left) and the standard gradient
norm (right) reach $10^{-3}, 10^{-4}$ and $10^{-5}$.
{\bf Top to bottom}: the rows show
optimization with a coarse mesh, a fine mesh, and the
same fine mesh uniformly scaled in dimension by $100\times$.
Note the consistency across mesh resolution and scaling
for our characteristic norm and
the disparity across the standard gradient norm.}
\label{fig:term_compare_2}
\end{figure}

Figure~\ref{fig:term_compare_2} illustrates consistency across
changing tolerance values, mesh resolutions, and scales.
example.  We show the iterates at measures
$10^{-3}$, $10^{-4} $and $10^{-5}$ for both our characteristic gradient
norm and the raw gradient norm, for meshes with varying refinement and
varying dimension (rescaling coordinates by a large factor). Similar
to Figure~\ref{fig:term_compare_1} comparing the vertex-normalized gradient
norm, there are large disparities for the raw gradient norm, but our
characteristic gradient norm is consistent.

\paragraph{Tolerances} The Swirl and Hilbert curve examples are
both extreme stress tests that require passing through low curvature
regions to transition from unfolding to folding; see e.g.,
Figure~\ref{fig:swirl_term_compare} above and our videos. For these
extreme tests we used a tolerance of $10^{-6}$ for our characteristic
gradient norm to consistently reach the final target shape.
However, for most practical geometry optimization tasks such a tolerance
is excessively precise. In experiments across a wide range of
energies and UV parametrization, 2D and 3D deformation tasks,
including those detailed below, we found that
$\|\nabla E(x)\| \leq 10^{-3} \langle W \rangle \| \ell(V,T) \|$
consistently obtained good-looking solutions with essentially no
further visible (or energy value) improvement possible. We
argue this is a sensible default except in pathological examples.
For all examples discussed here and below, with the exception of the
Swirl and the Hilbert curve tests, we thus use $\epsilon=10^{-3}$ for
testing termination.

\begin{figure*}[h!]
\centering
\includegraphics[width=1\linewidth]{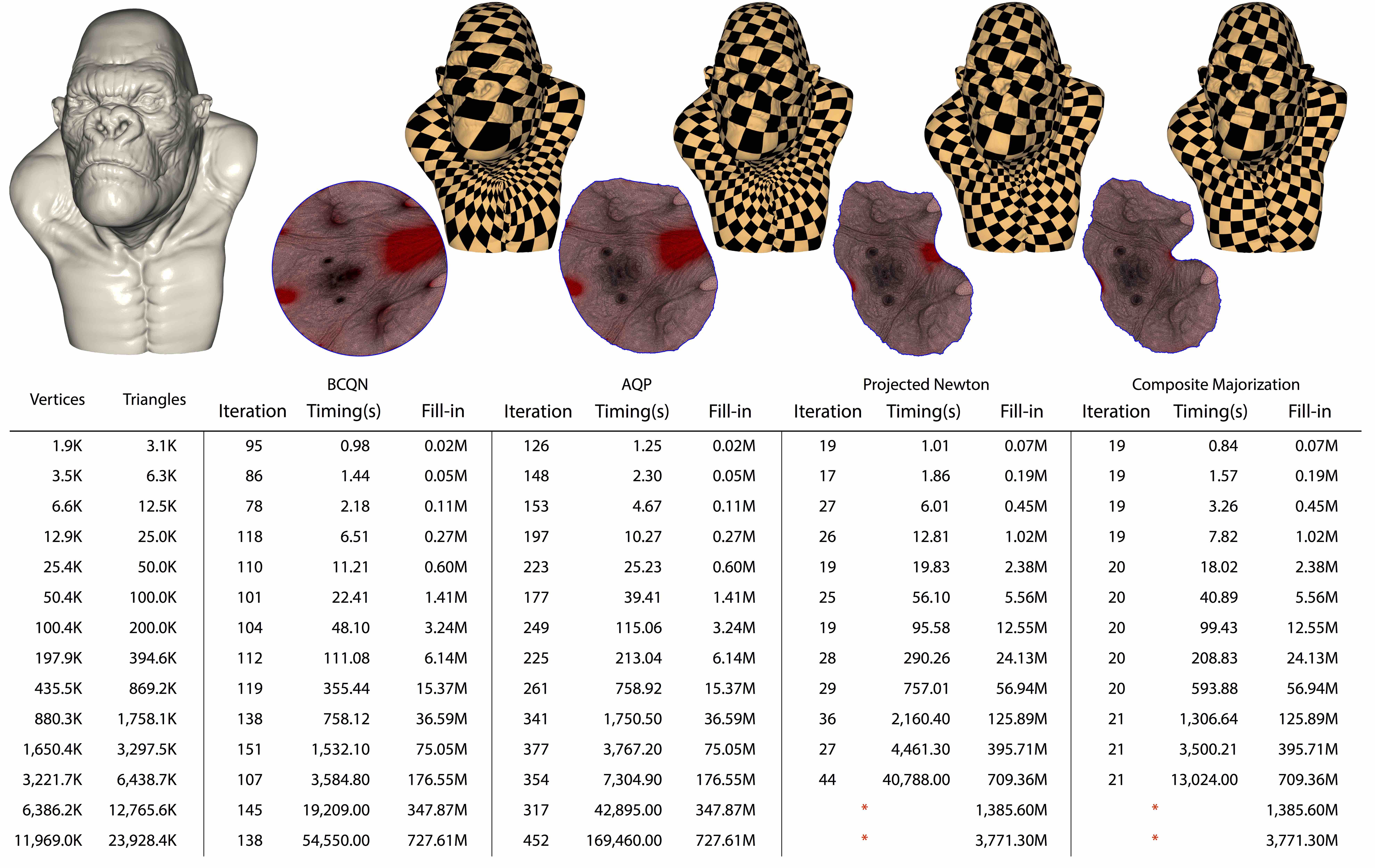}
\caption{\bfi{UV Parameterization Scaling, Timing and Sparsity.}
Performance statistics and memory use for increasing mesh
sizes up to 23.9M triangles, comparing BCQN with AQP, PN
and CM. For the Gorilla UV parametrization with ISO energy we
repeatedly double the mesh resolution and, for each method, report
number of iterations to convergence (characteristic norm $< 10^{-3}$),
wall-clock time (seconds) to convergence, and the nonzero fill-in
for the linear systems solved by each method. We use \red{\bf*} to
indicate out-of-memory failure for matrix factorization; see
\S\ref{sec:results} for discussion. Also note that stencils for CM
and PN are identical (differing only by actual entries) while AQP
and BCQN both solve with the same smaller scalar Laplacian.}
\label{fig:2d_scale_table}
\end{figure*}

\subsection{Newton-type methods}
\label{sec:newton-type}

\begin{figure}[h!]
\vspace{3mm}
\centering
\includegraphics[width=1\linewidth]{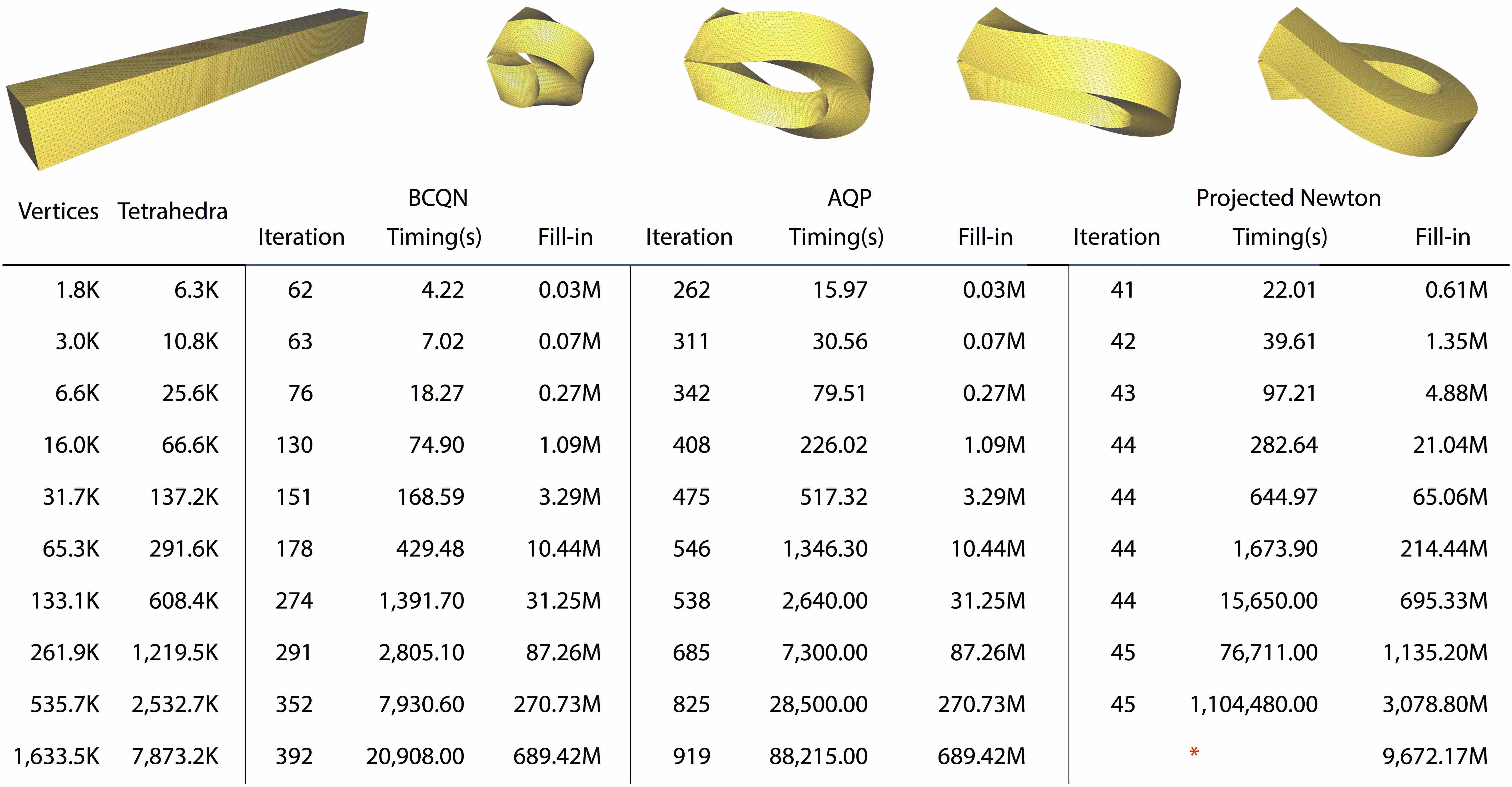}
\caption{\bfi{Three-Dimensional Deformation Scaling, Timing and
Sparsity.} Performance statistics and memory use for increasing
mesh sizes up to 7.8M tetrahedra, comparing BCQN with AQP
and PN. (CM does not extend to 3D.) We initialize a bar with a
straight rest shape to start in a tightly twisted shape, constraining
both ends to stay fixed and then optimize over increasing resolutions.
For each method we report number of iterations to convergence
(characteristic norm $< 10^{-3}$), wall-clock time (seconds) to
convergence, and the nonzero fill-in for the linear system solved
by each method. We use \red{\bf*} to indicate out of memory for the
computation on our test system; see \S\ref{sec:results} for discussion.
}
\label{fig:3d_scale_table}
\end{figure}

\begin{figure}[h!]
\centering
\includegraphics[width=1\linewidth]{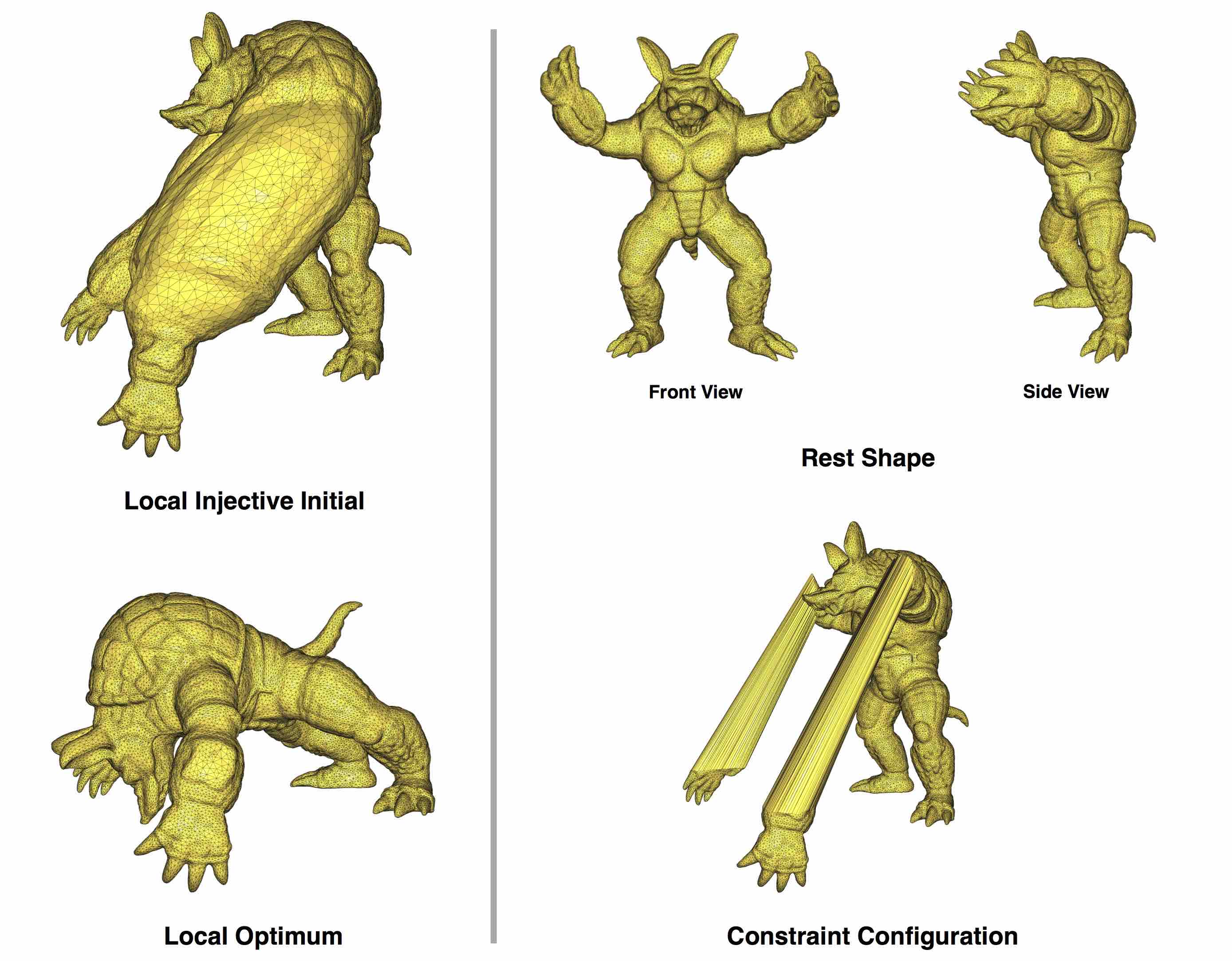}
\caption{\bfi{Armadillo Deformation Test}. We compare three-dimensional
deformation optimizations of a 1.5M element tetrahedral mesh of the
T-pose armadillo with BCQN and PN. We constrain the armadillo's
feet to rest position, its hands to touch the ground and use the
LBD method to create a locally injective
initialization for the solvers. Here BCQN requires 393 iterations
to converge while PN converges in just 9. However, as BCQN is much
cheaper and more scalable per iterate it takes only 4,148 seconds,
while PN spends 13,447 seconds.}
\label{fig:3d_large_defo}
\end{figure}

While Newton's method, on its own, handles convex energies like ARAP well~\cite{Chao:2010:ASG}
it is insufficient for nonconvex energies: modification of the Hessian is
required~\cite{Shtengel:2017:GOV,Nocedal:2006:Book}. Here we examine
the convergence, performance and scalability of Projected Newton (PN)~\cite{Teran:2005:RQF},
a general-purpose modification for nonconvex energies, and
CM~\cite{Shtengel:2017:GOV}, a more recent convex majorizer currently restricted
to 2D problems and a trio of energies (ISO, Symmetric ARAP and NH), and compare them with 
AQP and BCQN.
For the 2D parameterization problems in Figure~\ref{fig:2d_scale_table}
we can compare all four methods while for the 3D deformation
problems in Figures~\ref{fig:3d_scale_table} and \ref{fig:3d_large_defo}
CM is not applicable.

As we increase the size of the 2D problem by mesh refinement in
Figure~\ref{fig:2d_scale_table}, both CM and PN maintain low and almost constant
iteration counts to converge, with CM enjoying an advantage for larger problems;
in Figure~\ref{fig:3d_scale_table} 

Figures~\ref{fig:2d_scale_table}, \ref{fig:3d_scale_table}, and
\ref{fig:3d_large_defo} examine the scaling behavior of the various
methods under mesh refinement, for 2D parameterization and 3D deformation.
The Newton-type methods PN and CM (when applicable) maintain low
iteration counts that only grow slowly with increasing mesh size;
from the outset BCQN and AQP require more iterations, though the iteration count
also grows slowly for BCQN. Nonetheless, BCQN is the fastest across
all scales in each test as its overall cost per iteration remains much lower.
BCQN iterations require no re-factorizations (which scales poorly, particularly
in 3D, as discussed in Section \ref{sec:alg}) and only solves
smaller and sparser scalar Laplacian problems per coordinate compared
to the larger and denser system of CM and PN. This advantage for BCQN only
increases as problem size grows; indeed, for the largest problems BCQN
succeeded where CM and PN ran out of memory for factorization.

\subsection{A Note on Solving Proxies and Pardiso} 
\label{sec:pardiso}

Recent methods including CM have taken advantage of the efficiencies
and optimizations provided by the Pardiso solver.
While this can improve runtime of the factorization and backsolves
by a constant factor, it cannot change the asymptotic lower bounds on complexity;
the sparse matrix orderings in both SuiteSparse and Pardiso already appear
to achieve the bound on typical mesh problems.
In tests on our computer, across a large range of scales in two and
three dimensions, we found Pardiso was occasionally slower than
SuiteSparse but usually 1.4 to 3 times faster, and at most to 8.1 times faster
(for backsolving with a 3D scalar Laplacian).

Individual iterates of AQP have the same overall efficiency as BCQN (dominated
by the linear solves); switching to Pardiso leaves the relative performance of the two methods
unchanged. While CM and PN are even more dependent on the efficiency of the linear solver,
due to more costly refactorization each step, the same speed-ups possible with Pardiso also apply to
BCQN, so again there is no significant change in relative performance between the
methods.

\begin{figure}[h!]
\centering
\includegraphics[width=1\linewidth]{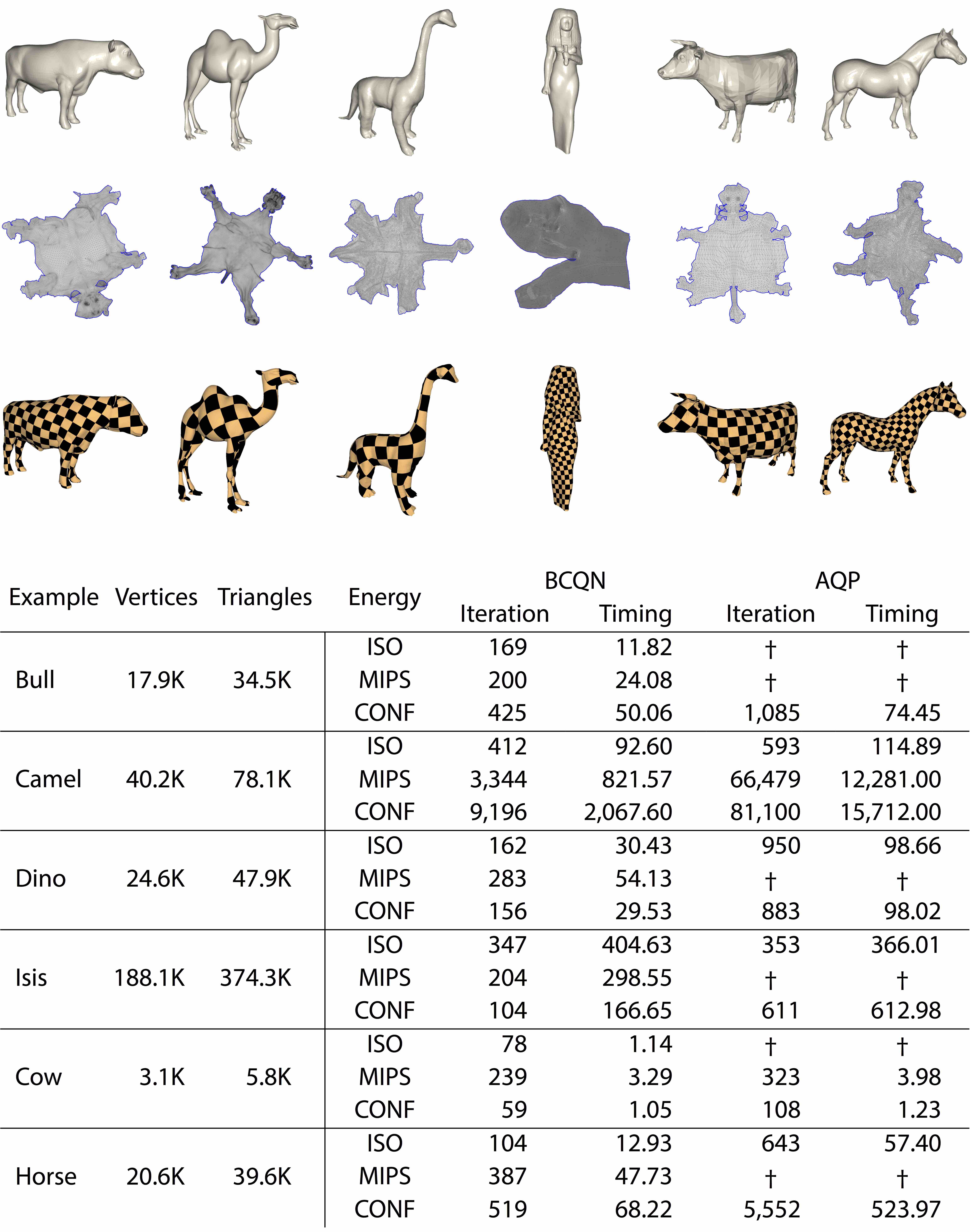}
\caption{\bfi{UV parameterization.} {\bf Top row:}  3D meshes  for
UV parametrization with ISO, MIPS, and CONF distortion energies.
{\bf Middle two rows:} converged maps and texturing from BCQN on
ISO examples. {\bf Bottom:} for each method / problem pair we report
number of iterations to convergence (characteristic norm $< 10^{-3}$)
and wall-clock time (seconds) to convergence. We use ${\bf \dagger}$
to indicate when AQP does not converge; see \S\ref{sec:1st}.}
\label{fig:uv_table}
\end{figure}

\begin{figure}[h!]
\centering
\includegraphics[width=1\linewidth]{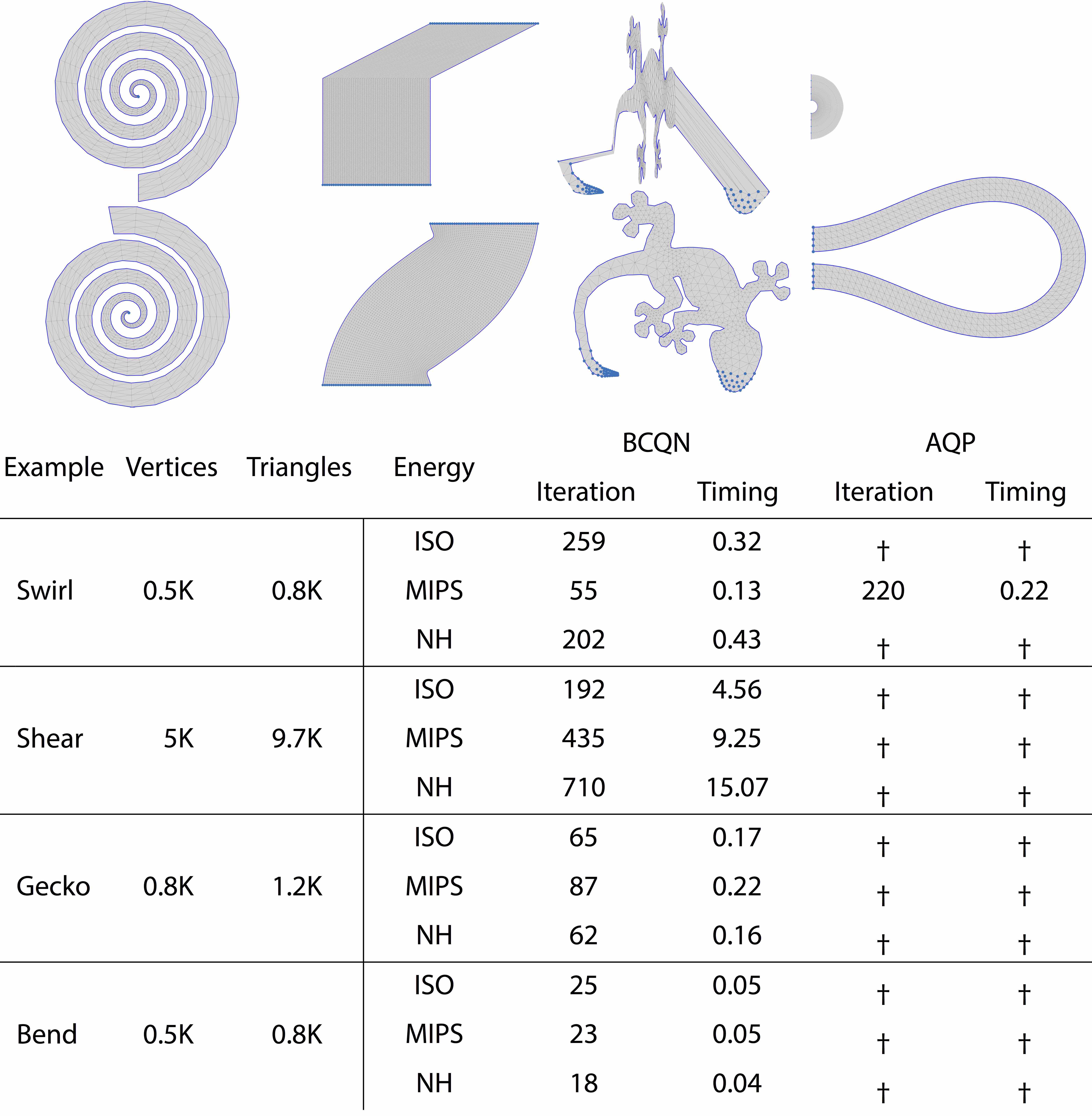}
\caption{\bfi{Two-Dimensional Deformation.} {\bf Top:}  initial
conditions and vertex constraints (blue points) for deformation
problems minimizing ISO, MIPS, and NH deformation energies. {\bf
Middle:} converged solutions from BCQN on ISO examples. {\bf Bottom:}
for each method / problem pair we report number of iterations to
convergence (characteristic norm $< 10^{-3}$) and wall-clock time
(seconds) to convergence. We use ${\bf \dagger}$ to indicate when
AQP does not converge; see \S\ref{sec:1st}.}
\label{fig:2d_defo_table}
\end{figure}

\begin{figure}[h!]
\centering
\includegraphics[width=1\linewidth]{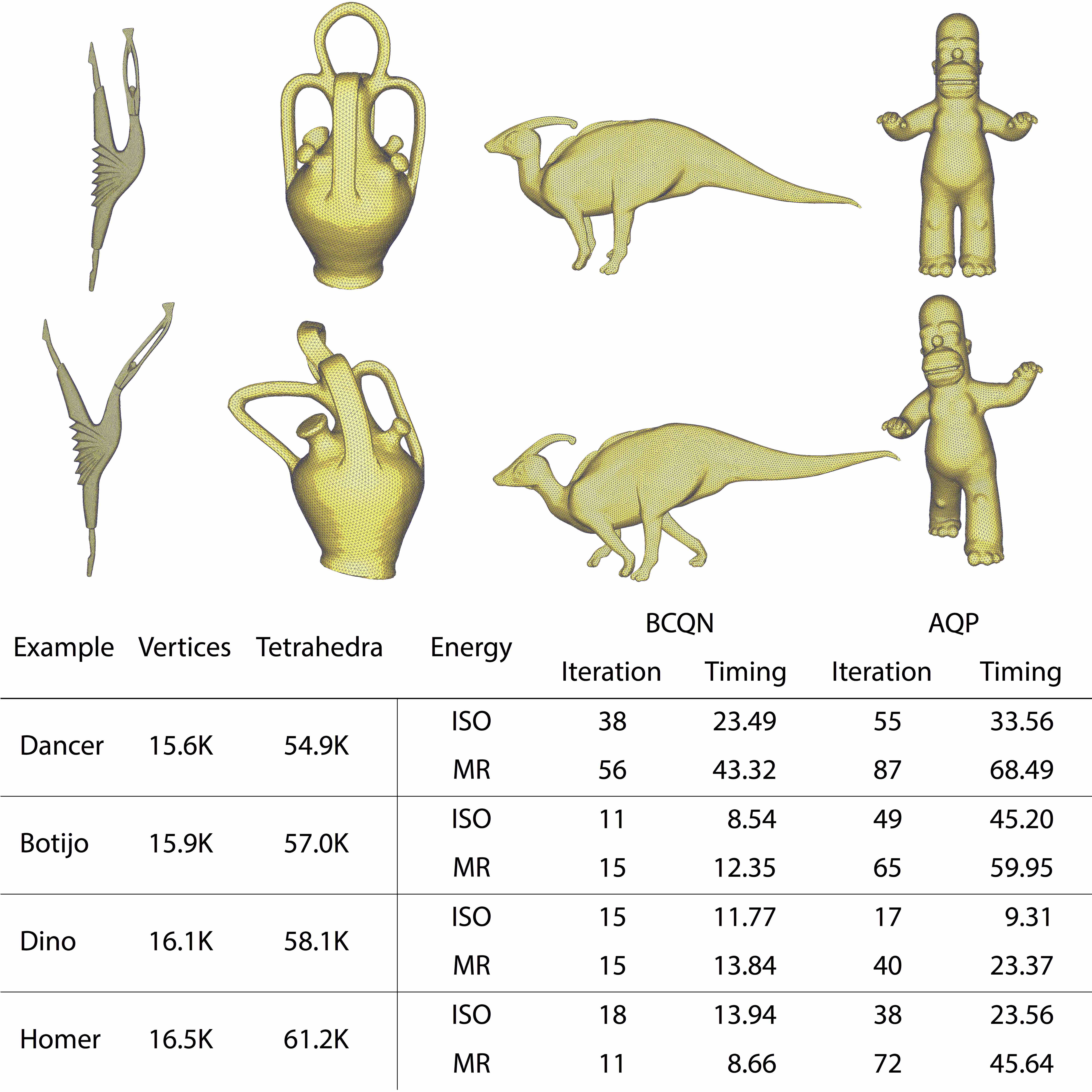}
\caption{ \bfi{Three-Dimensional Deformation.} {\bf Top:}  initial
conditions for vertex-constrained deformation problems minimizing
ISO and MR deformation energies. {\bf Middle:} converged solutions
satisfying constraints from BCQN on MR examples. {\bf Bottom:} for
each method / problem pair we report number of iterations to
convergence (characteristic norm $< 10^{-3}$) and wall-clock time
(seconds) to convergence. }
\label{fig:3d_defo_table}
\end{figure}


\subsection{First-order methods}
\label{sec:1st}

Among existing first-order methods for geometry optimization AQP
has so far shown best efficiency~\cite{Kovalsky:2016:AQP} with
improved convergence over SGD as well as standard L-BFGS.
Likewise, as we see in Figures\ \ref{fig:2d_scale_table},
\ref{fig:3d_scale_table}, and \ref{fig:3d_large_defo}, when we scale
to increasingly larger problems AQP will dominate over Newton-type
methods and so potentially offers the promise of reliability across
applications. Finally although small BQCN performs a small fixed
amount of extra work per-iteration in the line-search filter and
quasi-Newton update. Thus in Figures\ \ref{fig:2d_scale_table},
\ref{fig:3d_scale_table}, \ref{fig:uv_table} and \ref{fig:3d_defo_table}.
we compare AQP and BCQN over a range of practical geometry optimization
applications: respectively UV-parameterization, 2D deformation, and
3D deformation with nonconvex energies from geometry processing and
physics.  Throughout we note three key features distinguishing BCQN:

\bfi{Reliability and robustness.} 
AQP will fail to converge in some cases, see e.g. Figure\ \ref{fig:aqp_stop}, while BCQN reliably converges. In our testing AQP fails to converge in over 40\% of our tests with nonconvex energies; see e.g. Figures~\ref{fig:uv_table} and \ref{fig:2d_defo_table}.
This behavior is duplicated in our test-harness code and AQP's reference implementation.

\bfi{Convergence speed.} When AQP is able to converge, BCQN consistently provides faster convergence rates for nonconvex energies. In our experiments convergence rates range up to well over 10X with respect to AQP.

\bfi{Performance.} BCQN is efficient. When AQP is able to converge, BCQN remains fast with up to a well over 7X speedup over AQP on nonconvex energies.

\subsection{Across the Board Comparisons} 

Here we compare the performance and memory usage of BCQN with best-in-class geometry optimization methods across the board: AQP, PN and CM for both 2D parameterization and 3D deformation tasks. Results are summarized in Figures\ \ref{fig:2d_scale_table}, \ref{fig:3d_scale_table} and \ref{fig:3d_large_defo}. Note that CM does not extend to 3D.

In Figures~\ref{fig:2d_scale_table} and \ref{fig:3d_scale_table} we examine the scaling of AQP, PN, CM and BCQN to larger meshes and thus to larger problem sizes in both 2D parametrization (up to 23.9M triangles) and 3D deformation (up to 7.8M tetrahedra). As noted above: from the outset, BCQN requires more iterations than CM and PN; however, BCQN's overall low cost per iteration makes it faster in performance across problem sizes when compared to both CM and PN. We then note that AQP, on the other hand, has slower convergence and so, at smaller sizes it often does not compete with CM and PN. However, once we reach larger mesh problems, e.g. $\sim\geq$ 6M triangles in Figure~\ref{fig:2d_scale_table}, the cost of factorization and backsolve of the denser linear systems of CM and PN becomes significant so that even AQP's slower convergence results in improvement. This is the intended domain for which first-order methods are designed but here too, as we see in Figure~\ref{fig:2d_scale_table}, BCQN continues to outperform both AQP as well as CM and PN across all scales. Please see our supplemental video for visual comparisons of the relative progress of PN, CM, AQP and BCQN.

\section{Conclusion}

In this work we have taken new steps to both advance the state of
the art for optimizing challenging nonconvex deformation energies
and to better evaluate new and improved methods as they are
subsequently developed. Looking forward these minimization tasks
are likely to remain fundamental bottlenecks in practical codes and
so advancement here is critical. Our three primary contributions
together form the BCQN algorithm which pushes current limits in
deformation optimization forward in terms of speed, reliability, and
automatibility.
At the same time looking ahead we also expect that each contribution
individually should lead to even more improvements in the near
future.

\subsection{Limitations and Future Work}

While our focus is on recent challenging
nonconvex energies not addressed by the popular
local-global framework, similar to AQP we have observed
significant speedup for convex energies as well.  Currently
in comparing AQP and BCQN on the same set of 2D and 3D tasks with
the convex ARAP energy we observe a generally modest improvement
in convergence, up to a little over $4\times$, which is generally balanced
by the small additional overhead of BCQN iterations. Note for energies
like ARAP there is no barrier, hence no need for our line search filtering,
but other opportunities for improvement may be abailable in future research.

While our current blending model works well in
our extensive testing, it is empirically constructed; it is in no
sense proven optimal. We believe that further analysis, better
understanding and additional improvements in quasi-Newton blending
are all exciting and promising avenues of future investigation.

Finally, we note that while we have focused here on optimizing
deformation energies defined on meshes, there is a wide range
of critical optimization problems that take similar general structure:
minimizing separable nonlinear energies on graphs. Further extensions
are thus exciting directions of ongoing investigation.

\bibliographystyle{acmsiggraph}
\nocite{*}
\bibliography{paper}

\appendix

\section{Equivalence}
\label{sec:TH}

{\it Theorem.} For our energy densities $W(\sigma) = {f(\sigma)}/{g(\sigma)}$ with $f(\sigma)>0$ and $g(\sigma)\rightarrow 0$ as $\sigma\rightarrow 0$,
        $x^*$ is a stationary point of $\{E(x)   :   a(x) \geq 0 \}$ iff it is a locally injective stationary point of the unconstrained energy $E(x)$.

\begin{proof}
The $1/g(\sigma)$ term drives element energies $W\big(F_t(x)\big) \rightarrow \infty$ as $a_t(x) \rightarrow 0$. 
Stationary points $x_u^*$ of unconstrained $E$ are given by $\nabla E(x_u^*) = 0$ and must satisfy $|a(x_u^*)| > 0$. The addition of local injectivity then requires $a(x_u^*) > 0$.
Stationary points $x_c^*$ of $\{E(x) \>  : \>  a(x) \geq 0 \}$ are given by the Karush-Kuhn-Tucker (KKT) conditions
\begin{align}
\label{eq:nonlinear_KKT}
\nabla E(x_c^*) - \nabla a(x_c^*) \lambda = 0 \> \> \text{and} \> \> 
0 \leq \lambda \perp  a(x_c^*) \geq 0.
\end{align}
(Here $\lambda = (\lambda_1, .. ,\lambda_m)^T \in R^m$ is a Lagrange multiplier and $\vc x \perp \vc y$ is the \emph{complementarity condition} $y_t z_t = 0,\ \forall t$.)
All $x_u^*$ satisfy (\ref{eq:nonlinear_KKT}) with $\lambda > 0$. For $x_c^*$ satisfying (\ref{eq:nonlinear_KKT}) any $\lambda_t = 0 \implies a_t(x_c^*) = 0\implies W \big(F_t(x_c^*)\big)  = \infty$. Thus we must have  $\lambda > 0 \implies  a(x_c^*) > 0$ so that $x_c^*$ are locally injective stationary points of the unconstrained energy $E(x)$.
\end{proof}

\end{document}